\documentclass[final,3p,times,twocolumn]{elsarticle}

\usepackage{amsmath,amssymb,amsthm}
\usepackage{xcolor}
\usepackage{bm}
\usepackage[shortlabels]{enumitem}
\usepackage{lineno,hyperref} 
\usepackage{tikz-cd}
\modulolinenumbers[5]
\journal{Physica D}

\definecolor{armygreen}{rgb}{0.29, 0.33, 0.13}

\newcommand{\mA}{\mathcal{A}}
\newcommand{\mX}{\mathcal{X}}
\newcommand{\mY}{\mathcal{Y}}
\newcommand{\mE}{\mathcal{E}}
\newcommand{\mH}{\mathcal{H}}
\newcommand{\mG}{\mathcal{G}}

\newcommand{\BR}{\mathbb{R}}
\newcommand{\BE}{\mathbb{E}}
\newcommand{\BL}{\mathbb{L}}

\newcommand{\bk}{{\bf k}}

\newtheorem{thm}{Theorem}
\newtheorem{lem}[thm]{Lemma}
\newtheorem{prop}[thm]{Proposition}
\newdefinition{defn}[thm]{Definition}
\theoremstyle{remark}
\newtheorem{rk}[thm]{Remark}

\DeclareMathOperator{\spn}{span}
\DeclareMathOperator{\tr}{tr}
\DeclareMathOperator{\var}{var}
\DeclareMathOperator{\ran}{ran}
\DeclareMathOperator{\supp}{supp}
\DeclareMathOperator{\Leb}{Leb}
\DeclareMathOperator{\Id}{Id}

\begin{document}

%\linenumbers
\begin{frontmatter}

\title{Operator-Theoretic Framework for Forecasting Nonlinear Time Series with Kernel Analog Techniques}
\author[]{Romeo Alexander}
\author[]{Dimitrios Giannakis\corref{mycorrespondingauthor}}
\ead{dimitris@cims.nyu.edu}
\cortext[mycorrespondingauthor]{Corresponding author}
\address{Center for Atmosphere Ocean Science, Courant Institute of Mathematical Sciences, New York University, New York, New York 10012, USA}

\biboptions{sort&compress}
\bibliographystyle{elsarticle-num}

\begin{abstract}
    Kernel analog forecasting (KAF), alternatively known as kernel principal component regression, is a kernel method used for nonparametric statistical forecasting of dynamically generated time series data. This paper synthesizes descriptions of kernel methods and Koopman operator theory in order to provide a single consistent account of KAF. The framework presented here illuminates the property of the KAF method that, under measure-preserving and ergodic dynamics, it consistently approximates the conditional expectation of observables that are acted upon by the Koopman operator of the dynamical system and are conditioned on the observed data at forecast initialization. More precisely, KAF yields optimal predictions, in the sense of minimal root mean square error with respect to the invariant measure, in the asymptotic limit of large data. The presented framework facilitates, moreover, the analysis of generalization error and quantification of uncertainty. Extensions of KAF to the construction of conditional variance and conditional probability functions, as well as to non-symmetric kernels, are also shown. Illustrations of various aspects of KAF are provided with applications to simple examples, namely a periodic flow on the circle and the chaotic Lorenz 63 system. 
\end{abstract}


\end{frontmatter}

\section{Introduction}

Forecasting dynamically generated time series is a challenging problem that often requires statistical methods, especially when the underlying equations are either unknown or computationally intractable.
Data-driven forecasting methods have been sought after at least since Lorenz attempted to use naturally occurring historical analogs for climate predictions in the 1960's~\cite{Lorenz69b}. That early attempt was limited in success, but larger data sets and improved computing resources have made more recent analog-based nonparametric methods more viable~\cite{FarmerSidorovich87,casdagli1989nonlinear, SugiharaMay90, Sauer92, fan2008nonlinear, ZhaoGiannakis16}. Various types of ensemble analog forecasting are employed in short-term meteorological forecasts~\cite{delle2011kalman,AtenciaZawadzki15}, and versions of analog forecasting that utilize kernels have been shown to have predictive value for certain weather and climate phenomena~\cite{VanDenDool06,AlexanderEtAl17,ComeauEtAl17,DingEtAl18,ComeauEtAl19,WangEtAl19b}. 

While naturally occurring analogs may be a point of emphasis for nonparametric methods in physical science applications,  abstract statistical structures are the focus when situated in a more general machine learning context. Common nonparametric machine learning techniques include multilayer perceptrons \cite{voyant2015meteorological}, Bayesian neural networks~\cite{chakraborty1992forecasting}, classification and regression trees (CART), and a variety of kernel methods~\cite{ahmed2010empirical}. Although each of these methods can provide value in unique ways to specific problems, kernel methods are particularly well suited to problems where there may be a natural, a priori, notion of similarity between data points. Since analog methods rely on the possibility that the relevance of any historical analog to present day conditions can be quantitatively determined, formal understanding of such methods can improve when they are cast within the larger framework of kernel methods.

Kernel methods constitute a class of algorithms that perform classical calculations in a rich functional feature space in order to extract and predict nonlinear patterns. This central idea, commonly referred to as ``the kernel trick", was first proposed in 1964~\cite{aizerman1964theoretical}, was popularized with the invention of nonlinear support vector machines (SVMs) in 1992~\cite{boser1992training}, and has since spread to a variety of machine learning applications~\cite{hofmann2008kernel,SteinwartChristmann08}. Kernel methods for regression, such as support vector regression (SVR)~\cite{drucker1997support}, kernel ridge regression (KRR)~\cite{saunders1998ridge}, and kernel principal component regression (KPCR)~\cite{rosipal2000kernel}, may be applied to appropriately lagged signals to produce time series forecasts, such as with SVR forecasting~\cite{muller1997predicting}, KRR forecasting~\cite{exterkate2016nonlinear}, and KPCR forecasting~\cite{rosipal2000kernel}. Such kernel forecasting methods have been frequently used in finance and econometrics ~\cite{tay2001application}, and have recently found use in climate science~\cite{ZhaoGiannakis16, AlexanderEtAl17, ComeauEtAl17, ComeauEtAl19}, where they were termed kernel analog forecasting (KAF).

Statistical learning theory~\cite{vapnik1999overview} is the standard theoretical framework for deriving and analyzing kernel methods, among other machine learning algorithms. The learning guarantees and estimates of rates of convergence are well known when the underlying data are independently and identically distributed~\cite{cucker2007learning}. For time series, where the i.i.d.\ assumption is generally not valid, an extension of the standard i.i.d.\ statistical learning framework to that of stochastic processes has yielded softer guarantees that depend on mild conditions on the stationarity of the system~\cite{kuznetsov2015learning, kuznetsov2016time}. Although trajectories of a dynamical system can be viewed as a special case of a stochastic process, it is also worthwhile to employ the typical measure and operator-theoretic perspectives of modern dynamical systems theory~\cite{Baladi00,eisner2015operator}, where the induced action of the dynamical system on an intrinsically linear space of observables is given a more prominent role. This operator-theoretic perspective, although widespread in the study of dynamical systems \cite{BudisicEtAl12}, has yet to be fully exploited in conjunction with kernel forecasting methods.

The main contribution of this paper is a rigorous reformulation of KAF techniques within the framework of operator-theoretic ergodic theory and statistical learning theory. This view relies on the equivalence of forecasts with conditional expectation or, alternatively, geometric projection, both of which draw on the rich theory of functional analysis. One benefit from such a perspective is that it turns the problem of error analysis into the well studied problem of convergence in Hilbert spaces. Another benefit to this approach is that it demystifies the kernel functions somewhat by revealing their special role in bridging the gap between $L^2$ Hilbert spaces and the continuous function space in which forecasts are ultimately expressed. A third benefit is the modularity and extensibility that comes from casting kernel forecasting algorithms as a composition of operators applied to a careful choice of observable. In particular, by expressing forecasts as a composition of a regressor operator and the Koopman operator \cite{Koopman31}, the latter being a construct representing the action of evolving forward in time, features of the statistics and the dynamics are more easily separated and studied independently. For example, approximations of Koopman and the related transfer operators has been the subject of recent research~\cite{DellnitzJunge99,MezicBanaszuk04,Mezic05,RowleyEtAl09,Schmid10,BudisicEtAl12,FroylandEtAl13,FroylandEtAl14,GiannakisEtAl15,WilliamsEtAl15,KlusEtAl16,BruntonEtAl17,KlusEtAl18,KordaEtAl18,DasGiannakis19,Giannakis19}, and may be combined with approximations of the regressor operator to yield new formulations. Moreover, with appropriate choices of the response observable, forecasts can be obtained not just for the conditional mean of an observed quantity, but also that quantity's conditional variance and higher-order moments, which are important for uncertainty quantification. Conditional probability may also be approximated and predicted with a kernel analog approach. In this analysis, reproducing kernel Hilbert spaces (RKHSs) \cite{FerreiraMenegatto13,Paulsen16} play a central role as ambient hypothesis spaces of functions, with enough structure to enable an explicit representation of the forecasting function (also known as target function) in a fully empirical manner.

This paper is organized as follows. Section~\ref{secKAF} introduces the forecasting problem under study, and describes the KAF framework. Section~\ref{secError} studies the generalization error of constructed forecasts, paying particular attention on how to quantify the discrepancy between empirical and ideal forecasts. Our main result on the convergence of KAF to the conditional expectation is stated as Theorem~\ref{thmEmpiricalConv} in that section. Section~\ref{secExt} introduces a few extensions, including KRR, non-symmetric kernels, conditional variance, and conditional probability. In Section~\ref{secKernel}, we provide general guidelines for choosing the kernel. Section~\ref{secApps} shows the result of applying KAF to two examples, namely a periodic flow on the circle and the chaotic Lorenz 63 (L63) system \cite{Lorenz63}. Section~\ref{secConclusions} provides our principal conclusory remarks, and examines the applicability of KAF to various real-world problems. Technical results are collected in \ref{app}.

\section{\label{secKAF}Kernel analog forecasting (KAF) techniques}

In this section, we describe the mathematical framework underlying the KAF approach introduced in \cite{ZhaoGiannakis16}. We start from a general formulation of forecasting as error minimization (Section~\ref{secBackground}), and gradually build onto that various dynamical systems and functional analytic tools, leading (in Section~\ref{secTarget}) to the construction of the RKHS-based KAF target function. It should be noted that our exposition differs substantially from \cite{ZhaoGiannakis16}, which focuses heavily on RKHS interpolation theory from the outset. In particular, an advantage of the perspective put forward here is that the RKHS formalism emerges as a natural consequence of seeking target functions in an explicitly constructible ambient hypothesis space with a Hilbert space structure, as opposed to the more ``axiomatic'' use of RKHSs in \cite{ZhaoGiannakis16}. This perspective will also facilitate the error analysis in Section~\ref{secError}.

\begin{figure*}
    \centering\includegraphics[width=\linewidth]{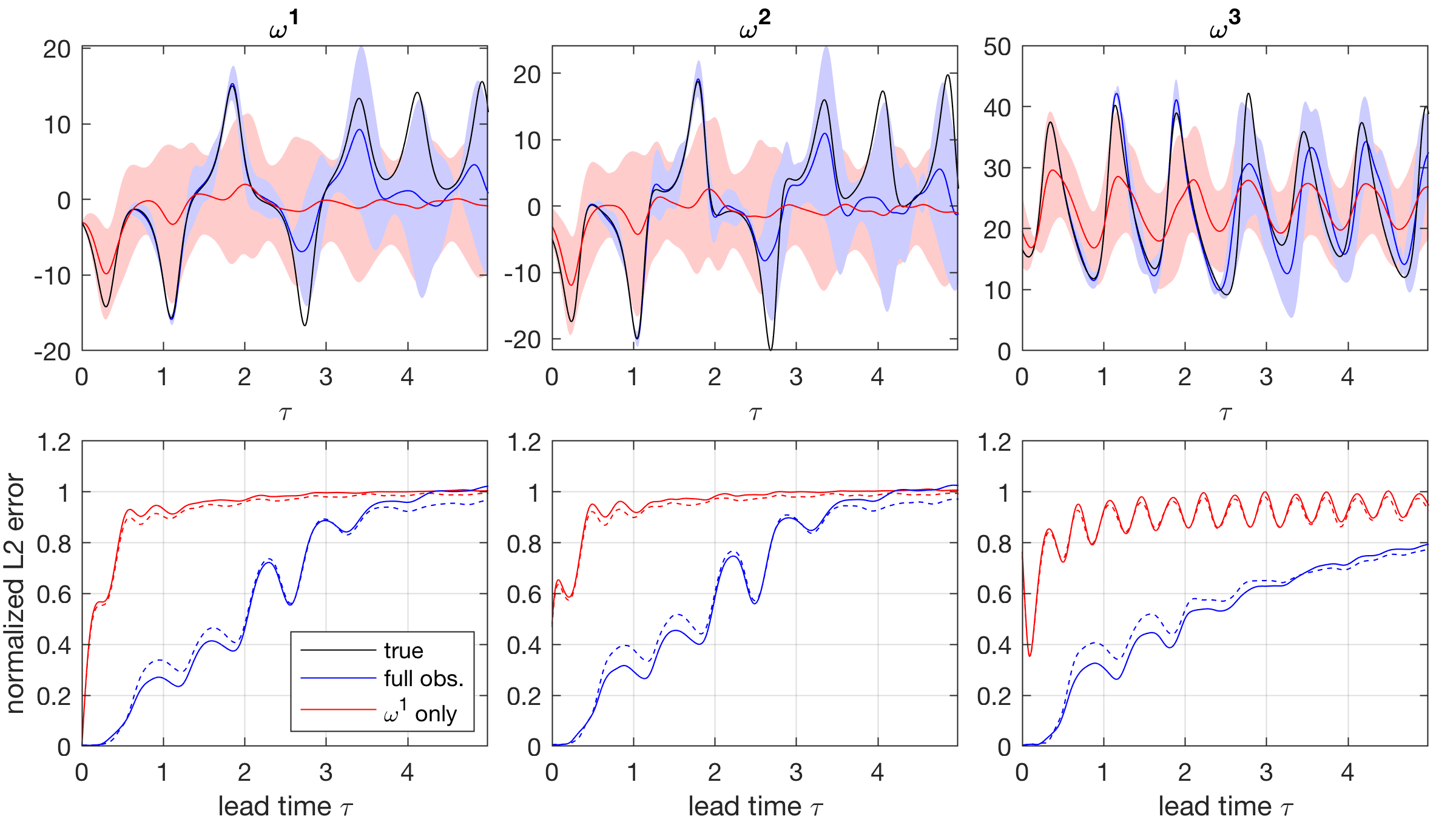}
    \caption{\label{figL63Pred}Results of KAF applied to prediction of the components of the state vector $ \omega = (\omega^1,\omega^2,\omega^3) \in \mathbb R^3$ of the L63 system, under full (blue lines) and partial (red lines) observations. In the fully observed case, the covariate $X$ is the identity map on $\Omega = \mathbb R^3$. In the partially observed case, $X(\omega) = \omega^1$ is the projection to the first coordinate. Both cases utilize training datasets of $n=\text{64,000}$ samples $(x_j,y_j)$, with $ x_j = X(\omega_j) $, $y_j = Y(\omega_j)$, obtained from the same underlying L63 trajectory $\omega_j = \Phi^{(j-1)\,\Delta t}(\omega_1) \in \Omega $ at a sampling interval of $\Delta t = 0.01$ natural time units. Top panels: True evolution $U^\tau Y(\omega)$ of the response $Y(\omega) = \omega^j$ (black lines) and forecast trajectories $f_{\tau,\ell,n}(x)$ (solid blue and red lines) as a function of lead time $ \tau $, obtained via the empirical target functions in Definition~\ref{defnEmpTarget} with $ \ell = 3000$ (full observations) or $\ell = 1000$ (partial observations). The forecasts are initialized from a fixed initial condition $ x = X(\omega) $ in the verification dataset. Shaded regions show error bounds obtained by adding $\pm \varepsilon_{\tau,\ell,n}(x)$ to the forecast trajectories, where $\varepsilon_{\tau,\ell,n}(x)$ is a KAF-derived estimate of the conditional standard deviation, given by~\eqref{eqErrEst} using the same training parameters as the $Y$ forecasts. Bottom panels: Root mean square error (RMSE) as a function of $ \tau $ determined from a verification dataset of $\tilde n = \text{64,000}$ (solid lines). The RMSE is normalized by the empirical standard deviation of $ \omega^j $ on the training dataset. Dashed lines show the normalized RMS value of the estimated error $ \varepsilon_{\tau,\ell,n} $. The agreement between actual and estimated errors indicates that $\varepsilon_{\tau,\ell,n}$ provides useful uncertainty quantification.}
\end{figure*}

Figure~\ref{figCommut} depicts the relationships between the function spaces and operators involved in the construction of the KAF target function in the form of a commutative diagram. The basic steps of the construction are also summarized in pseudocode in Table~\ref{tablePseudocode}. Figure~\ref{figL63Pred} shows an application of KAF to the L63 system under full and partial observations. This L63 application provides a guiding example of a number of challenges encountered in statistical forecasting, including partial state observations, mixing (i.e., chaotic) dynamics, and invariant measures supported on non-smooth attractors.   

\begin{figure*}
    \begin{displaymath}
        \begin{tikzcd}
            &&& C(\Omega) \arrow[dll, "\iota", blue] \arrow[drr, "\iota_n"', red] \arrow[d, "U^\tau"'] && & \\
 & L^2(\mu) \arrow[d, "U^\tau",blue] && C(\Omega) \arrow[lld, "\iota"'] \arrow[dd,"\iota"']&& L^2(\mu_n) \arrow[d, "U^{q}_n", red]&  \\
 &L^2(\mu)\arrow[d, "\Pi_X", blue] &  &  &   & L^2(\mu_n) \arrow[d, "\Pi_{X,n}",red]& \\
 &L^2_X (\mu) \arrow[ddl, shift left=1ex, armygreen, "\Xi^*"] \arrow[rr,dashed,blue] && L^2(\mu) & & L^2_X (\mu_n )\arrow[ddr, shift right=1ex, red, "\Xi^*_n"']&  \\
  & & L^2_X(\mu) \arrow[ul, dashed, armygreen] \arrow[ur,dashed,armygreen]  &  &  L^2_X(\mu) \arrow[ll, red, dashed] \arrow[ul,red,dashed]& &  \\
 L^2(\mu_X) \arrow[uur, shift left=1ex, "\Xi"] \arrow[rd, "T_\ell", armygreen]&  &  L^2(\mu_X) \arrow[u, armygreen,"\Xi"] & & L^2(\mu_X) \arrow[u, red, "\Xi"]  &  &  L^2(\mu_{X,n})\arrow[uul, shift right=1ex,"\Xi_n"'] \arrow[dl,"T_{\ell,n}", red]\\
 & \mathcal H_\ell  \arrow[r, armygreen,"\subset"]& C(\mathcal X) \arrow[u, armygreen, "\iota"] & & C(\mathcal X)\arrow[u, red, "\iota"']\arrow[ll, red, dashed] & \mathcal H_{\ell, n}\arrow[l,  red,"\subset"']& 
        \end{tikzcd}
    \end{displaymath}
    \caption{\label{figCommut} Commutative diagram illustrating the relationships between the function spaces and operators employed in KAF. We start from the space of continuous functions on $\Omega$, wherein the response variable $Y$ lies. The path proceeding downwards from $C(\Omega)$ along the black arrows represents the evolution of continuous functions under the Koopman operator $U^\tau $ for lead time $\tau$, followed by inclusion in the $L^2(\mu)$ space associated with the invariant measure. This leads to the observable $ U^\tau Y \in L^2(\mu)$, which we seek to approximate with minimal error with respect to $L^2(\mu)$ norm. The paths demarcated by blue, green, and red arrows represent three levels of approximation, each with its own errors (indicated by dashed arrows). The path colored in blue leads to the conditional expectation $\mathbb E[ U^\tau Y \mid X] = \Pi_X U^\tau Y $, given by the orthogonal projection of $U^\tau $ into the Hilbert subspace $L^2_X(\mu) \subseteq L^2(\mu)$ consisting of pullbacks of functions on covariate space $\mathcal X$. The conditional expectation is associated with the regression function $Z_\tau = \Xi^* \mathbb E[U^\tau Y \mid X]$ from Definition~\ref{defnRegression}, and can exhibit an irreducible form of error relative to $U^\tau Y $ if $X$ is not injective (indicated by a dashed blue arrow). The green-colored path describes the approximation of the regression function by a continuous function $ f_{\tau,\ell} $ on covariate space (the ideal target function in Definition~\ref{defnTarget}), lying in a finite-dimensional hypothesis space $\mathcal H_\ell$, which is a subspace of an RKHS $\mathcal H$ on $\mathcal X$. The operator $T_\ell$ leading to the target function is a spectrally truncated Nystr\"om extension operator, where $\ell$ is the number of kernel eigenfunctions  employed. This imparts an additional error, which vanishes, however, as $\ell \to \infty$ if the reproducing kernel of $\mathcal H$ is strictly positive definite. The path colored in red represents a data-driven approximation $f_{\tau,\ell,n}$ of $f_{\tau,\ell}$ (the empirical target function in Definition~\ref{defnEmpTarget}), obtained by replacing the invariant measure $ \mu $ by the sampling measure $ \mu_n $ associated with a training dataset consisting of $n$ samples, and the Koopman operator $U^{\tau} $ by the shift operator $U^{q}_n $ on the $n$-dimensional Hilbert space $L^2(\mu_n)$. Here, $ \tau = q \, \Delta t$, where $ \Delta t $ and $q$ are the sampling interval and number of timesteps, respectively. Note that there is no path connecting the Koopman operator on $C(\Omega)$ with the shift operator on $L^2(\mu_n)$ as in the $L^2(\mu)$ case since the dynamical flow is singular with respect to the sampling measure. The estimator $f_{\tau,\ell,n}$ exhibits a sampling error relative to $ f_{\tau,\ell} $ (indicated by dashed red arrows) which vanishes almost surely as $ n \to \infty $ by ergodicity.}
\end{figure*}

\begin{table*}
    \caption{\label{tablePseudocode}Pseudocode outlining the construction of the empirical KAF target function $f_{\tau,\ell,n}$ from Definition~\ref{defnEmpTarget}.}
    \small
    \hrulefill
    \begin{itemize}
    \item Inputs 
    \begin{itemize}
        \item Covariate training data $ x_1, \ldots, x_{n} \in \mathcal X $ at sampling interval $\Delta t$ 
        \item Response training data $ y_1, \ldots, y_{n} \in \mathcal Y $ at sampling interval $\Delta t$
        \item Symmetric positive-definite kernel function $ k_n : \mathcal X \times \mathcal X \to \mathbb R$
        \item Forecast timesteps $ q \in \mathbb N_0$
        \item Number of principal components (eigenfunctions) $ \ell \in \mathbb N$
    \end{itemize}
    \item Outputs 
        \begin{itemize}
            \item Target function $f_{\tau,\ell,n} : \mathcal X \to \mathbb R $ for lead time $ \tau = q \, \Delta t$ 
        \end{itemize}
    \item Steps 
    \begin{enumerate}
        \item Compute the leading $\ell $ eigenvectors $ \bm \phi_j \in \mathbb R^N $ of the $n \times n $ kernel matrix $ \bm K = [ k_n(x_i, x_j) ]$, arranged in order of decreasing corresponding eigenvalue $\lambda_{j,n}$. Normalize the eigenvectors such that $ \bm \phi_i \cdot \bm \phi_j / n = \delta_{ij} $.
        \item Form the $q$-step shifted response vector $ \bm y_{\tau} = ( 0, \ldots, 0, y_1, \ldots, y_{n-q} ) \in \mathbb R^n $, and compute the expansion coefficients $ \alpha_j(\tau) = \bm \phi_j \cdot \bm y_{\tau} / n $ for $ j \in \{ 1, \ldots, \ell \}$.  
        \item Form the orthonormal RKHS functions $ \psi_{j,n}(x) = \bm k(x) \cdot \bm \phi_j / ( n \lambda_{j,n}^{1/2} ) $, where $ \bm k(x) = ( k_n( x, x_1 ), \ldots, k_n(x,x_n) ) $ is the kernel vector.
        \item Form the target function $ f_{\tau,\ell,n}(x) = \sum_{j=1}^{\ell} \alpha_{j,n}(\tau) \psi_{j,n}(x) $.
        \end{enumerate}
    \end{itemize}
    \hrulefill
\end{table*}

\subsection{\label{secBackground}Mathematical background}
\paragraph*{Measure-theoretic framework}
In the measure-theoretic setup that we wish to pursue here, the primary object is a probability space $(\Omega, \mathcal F, \mu)$, where $\Omega$ is the space of all possible initial states, $\mathcal F$ is a $\sigma$-algebra of distinguished subsets of $\Omega$, and $\mu : \mathcal F \to \BR$ is a probability measure. We also have a measurable covariate space $(\mX, \Sigma_\mX)$, a measurable response space $(\mY, \Sigma_\mY)$, and, for each time $t \geq 0$, data-producing measurable functions $X_t : \Omega \to \mX$ and $Y_t : \Omega \to \mY$. By data-producing, we mean that the covariate and response data, $x_t$ and $y_t$, are regarded as the output of $X_t$ and $Y_t$, respectively, so that $x_t = X_t(\omega)$ and $y_t = Y_t(\omega)$ for some $\omega \in \Omega$. The space $\mY$ is assumed to be a Hilbert space over the complex numbers, whose inner product, $\langle \cdot, \cdot \rangle_\mY$, is taken to be conjugate-linear in its first argument. Note that we do not require that the space $\mX$ be linear.

The task of forecasting is to produce a measurable function $f_\tau : \mX \to \mY$ for any given lead time $\tau \ge 0$, referred to as the target function, such that $f_\tau \circ X_t$ approximates $Y_{t + \tau}$. A heuristic for selecting such an approximation is the variational approach, wherein $f_\tau$ is viewed as a minimizer of some global measure of error. The mean square error is a common such functional, given, as we will see below, its connection to Hilbert space theory. In particular, we regard $Y_t$ as an element of the space $\BL^2(\mu)$ of functions  $\Omega \to \mY$ that are square-integrable with respect to $\mu$. The target function $f_\tau$, meanwhile, is sought in the space $\BL^2(\mu_{X_t})$ of functions $f : \mX \to \mY$ that are square-integrable with respect to $\mu_{X_t}$, where $\mu_{X_t}$ is the pushforward of $\mu$ along $X_t$ (i.e., $\mu_{X_t}(S) = \mu(X_t^{-1}(S))$ for all $S \in \Sigma_X$). This implies that $ f\circ X_t $ is a square-integrable function in $\BL^2(\mu)$. 

In what follows, $L^2(\mu)$ will denote the Hilbert space of equivalence classes of functions in $\BL^2(\mu)$ taking $\mu$-a.e.\ equal values, equipped with the standard  inner product $\langle g_1, g_2 \rangle_{L^2(\mu)} = \int_\Omega \langle g_1(\omega), g_2(\omega) \rangle_\mY \, d \mu(\omega) $. We define the $L^2(\mu_{X_t})$ Hilbert spaces associated with $\BL^2(\mu_{X_t})$ analogously. As is customary, we will oftentimes identify functions in $\BL^2$ with their corresponding $L^2$ equivalence classes, but for the purpose of constructing concrete target functions we will keep elements of these spaces distinct. The mean square error of the target function $f_\tau$, given a lead time $\tau \ge 0$, may then be defined as
\begin{multline*}
    \lVert f_\tau \circ X_t - Y_{t+\tau} \rVert_{L^2(\mu)}^2 \\ 
    = \int_\Omega \lVert f_\tau \circ X_t( \omega)  - Y_{t + \tau}( \omega ) \rVert_\mY^2 \, d\mu(\omega). 
\end{multline*}

\paragraph*{Dynamical system framework}

A dynamical system on the space $\Omega$ is represented by a semigroup of measurable maps, $\{\Phi^t : \Omega \to \Omega \}_{t \ge 0}$, which evolve an initial state $\omega_0$ to a new state $\omega_t$. The function $X_t$ may then be represented by $X \circ \Phi^t$, where $X: \Omega \to \mX$. The response function $Y_{t + \tau}$ can be similarly broken up, but with the added step of using the flow map semigroup properties to split up $\Phi^{t + \tau}$ into $\Phi^\tau\circ \Phi^t$, resulting in the expression $Y_{t + \tau} = (Y \circ \Phi^\tau) \circ \Phi^t$, where $Y : \Omega \to \mY$. It is frequently useful to express the composition $Y\circ \Phi^\tau$ as the act of applying an operator $U^\tau$, known as the Koopman operator \cite{koopman1931hamiltonian}, on measurable $\mY$-valued functions on $\Omega$, so that $U^\tau Y = Y\circ \Phi^\tau$.  Note that, unlike $\Phi^\tau$, $U^\tau$ is an intrinsically linear operator.

Henceforth, we will assume that the dynamical system is measure-preserving; that is, the pushforward measure of $\mu$ along $\Phi^t$, denoted by $\mu_{t}$, is constant with respect to time, so that we may write $\mu_{t} = \mu$ for all times $t \geq 0$.  With such an assumption, the Koopman operator on measurable functions lifts to a unitary operator on $L^2(\mu) $, which we will denote using the same symbol $ U^\tau $. Moreover, the   mean square error is independent of the initialization time $t$, and is expressed as
\begin{equation}
\lVert f_\tau \circ X - U^\tau Y \rVert^2_{L^2(\mu)}.
\label{eq:MSE}
\end{equation}

\paragraph*{Conditional expectation}

The random variable $X$ induces a sub-$\sigma$ algebra $\mG \subseteq \mathcal F$, defined by $\mG = X^{-1} (\Sigma_{\mathcal X}) $. This means that every function $g : \Omega \to \mY$ that is measurable with respect to $\mG$ is such that $g = f \circ X$ for some $f : \mX \to \mY$. Thus, $\mathcal{G}$-measurable functions can be thought of as being ``coarser'' than $\mathcal{F}$-measurable functions, in the sense that they necessarily take constant values on subsets of $\Omega$ where $X$ is constant. We will denote the Hilbert subspace of $L^2(\mu)$ consisting of $\mathcal{G}$-measurable equivalence classes of functions by $L^2_X(\mu)$. The composition map by $X$, i.e., $ f \mapsto f \circ X $, then describes an isometric embedding $ \Xi : L^2(\mu_X) \to L^2(\mu)$, with range $L^2_X(\mu)$. 
It is a consequence of the Radon-Nikodym theorem \cite[][Chapter~5]{Kallenberg97} that for $U^\tau Y \in L^2(\mu)$,  there exists a unique $\mG$-measurable element $Z_\tau \circ X \in L^2(\mu)$, such that for all $g \in L^2_X(\mu)$,
\begin{equation}
\langle g, U^\tau Y\rangle_{L^2(\mu)} = \langle g, Z_\tau \circ X \rangle_{L^2(\mu)}.
 \label{eq:conditional_identity}
\end{equation}
It follows from this property that  $Z_\tau \circ X$ is the unique element in $L^2(\mu)$, or, equivalently, that $Z_\tau$ the unique element in $L^2(\mu_X)$, that minimizes the mean square error in~\eqref{eq:MSE}. We shall refer to the composition $Z_\tau \circ X$ as the conditional expectation $\BE[U^\tau Y \mid X]$ and to $Z_\tau$ as the regression function. It follows from the Hilbert space projection theorem that $ \mathbb{E}[U^\tau Y \mid X]$ has the geometrical interpretation of being the orthogonal projection of $U^\tau Y $ onto  $L^2_X(\mu) $.  That is, 
\begin{equation}
    \label{eqEProj}
    \mathbb{E}[U^\tau Y \mid X] = \Pi_X U^\tau Y,
\end{equation}
where $\Pi_{X} : L^2(\mu) \to L^2(\mu)$ is the orthogonal projection mapping into $L^2_X(\mu)$. Because the conditional expectation lies in $L^2_X(\mu)$, there exists a unique observable $Z_\tau \in L^2(\mu_X))$ on covariate space such that 
\begin{displaymath}
    \mathbb E[U^\tau Y \mid X ] = \Xi Z_\tau = Z_\tau \circ X. 
\end{displaymath}
This leads to the notion of the regression function, defined below through the adjoint map $ \Xi^* : L^2(\mu) \to L^2(\mu_X) $ with $\ker \Xi^* = ( L^2_X(\mu) )^\perp$. 

\begin{defn}[regression function]\label{defnRegression}The regression function at lead time $\tau$ associated with the response $Y$ and covariate $X$ is the $L^2(\mu_X)$ observable
\begin{displaymath}
    Z_\tau = \Xi^* U^\tau Y = \Xi^* \mathbb E[ U^\tau Y \mid X ]. 
\end{displaymath}
\end{defn}

By virtue of its error-minimizing properties, it is natural to seek forecasting algorithms producing target functions that consistently approximate $Z_\tau$. In the ensuing sections, we will show that under suitable ergodicity assumptions, KAF naturally produces such consistent estimators of the regression function from time-ordered samples of $X$ and $Y$ along a dynamical trajectory, without requiring prior knowledge of the underlying equations of motion.

\subsection{\label{secHypothesis}Hypothesis spaces}
\paragraph*{Learning framework}
Constructing the target function requires distinguishing between the spaces $\BL^2(\mu_X)$ and $L^2(\mu_X)$, which we do by way of the linear map $\iota : \BL^2(\mu_X) \to L^2(\mu_X)$ that associates each concrete function $f$ to its equivalence class $\iota f$. The mean square error is then represented with the functional $\mE_\tau : \BL^2(\mu_X) \to \BR$, known as the generalization error in machine learning contexts \cite{cucker2007learning}, defined by
\begin{equation}
    \label{eqGen}
    \mE_\tau(f) := \lVert \iota f \circ X - U^\tau Y\rVert_{L^2(\mu)}^{2}.
\end{equation}
The Hilbert space structure of $L^2(\mu)$, as well as the error-minimizing property of the conditional expectation $Z_\tau \circ X$, allows the generalization error to be decomposed as
\begin{displaymath}
\mE_\tau(f) = \mA_\tau(f) + \sigma_\tau,
\end{displaymath}
where $\mA_\tau(f)$ is the excess generalization error,
\begin{equation}
    \label{eqExcessGen}
\mA_\tau(f) = \lVert \iota f - Z_\tau \rVert^2_{L^2(\mu_X)},
\end{equation}
and $\sigma_\tau$ is the error intrinsic to the system and choice of covariate and response functions,
\begin{equation}
\label{eqSigmaRho}
    \sigma_\tau = \lVert Z_\tau \circ X - U^\tau Y \rVert^2_{L^2(\mu)}.
\end{equation}
Since $\sigma_\tau$ does not depend on $f$, minimizing $\mE_\tau$ is equivalent to minimizing $\mA_\tau$.

\paragraph*{Hypothesis space}
Constraints on the search for a minimizer of $\mA_\tau$ are characterized in terms of a hypothesis space $\mH \subseteq \BL^2(\mu_X)$ of functions. When the image $ H := \iota \mH $ is a closed and convex subset of the Hilbert space $L^2(\mu_X)$, then there exists a unique $g \in H$ such that $\inf_{h\in H} \lVert h - Z_\tau\rVert_{L^2(\mu_X)} = \lVert g - Z_\tau\rVert_{L^2(\mu_X)}$. Consequently, there exists  $f \in \mH$ for which $\iota f = g$, and thus $\inf_{f\in \mH} \lVert \iota f - Z_\tau\rVert_{L^2(\mu_X)} = \lVert g - Z_\tau\rVert_{L^2(\mu_X)}$. A sufficient condition for uniqueness of $f$ is that $\iota : \mH \to L^2(\mu_X)$ be an injection. 

\paragraph*{The pseudoinverse}
Assuming that $H$ is closed and convex in $L^2(\mu_X)$, so that there exists a well-defined orthogonal projection map $\Pi_{H} : L^2(\mu_X) \to L^2(\mu_X)$ mapping into $H$, the excess generalization error may be decomposed as
\begin{displaymath}
    \mA_\tau(f) = \left\lVert \iota f - \Pi_{H} Z_\tau \right\rVert^2_{L^2(\mu_X)} + \left\lVert \left(I - \Pi_{H} \right)Z_\tau \right\rVert^2_{L^2(\mu_X)}.
\end{displaymath}
The minimizer of $\mA_\tau$ over the hypothesis space $\mH$, therefore, is found by minimizing the norm of $\iota f - \Pi_{H} Z_\tau$. When $\iota$ is injective on $\mH$, then the restriction $\iota|_\mH$ of $\iota$ onto $\mH$ is invertible as a map $\iota|_\mH : \mH \to H$. In such a case, the unique minimizer of $\mA_\tau$ in $\mathcal{H}$ is expressible as
\begin{equation}
    \label{eqFTarget}
    f_{\tau,\mH} = \left(\iota|_\mH \right)^{-1} \Pi_{H} Z_\tau,
\end{equation}
and satisfies
\begin{equation}
    \mA_\tau(f_{\tau,\mathcal{H}}) = \lVert (I - \Pi_H) Z_\tau \rVert^2_{L^2(\mu_X)}.
    \label{eqExcessErrFl}
\end{equation}
We refer to this minimizer as the ideal target function since it does not depend on any training data.

\begin{defn}[ideal target function] \label{defnTarget} The ideal target function $f_{\tau,\mathcal H}$ at lead time $\tau$ associated with the response $Y$ and hypothesis space $\mathcal H$ is the minimizer of the excess generalization error functional $\mathcal A_\tau$ over $\mathcal H$, given by~\eqref{eqFTarget}.
\end{defn}

We shall refer to the map $T : L^2(\mu_X) \to \mH $, with $ T = \left(\iota|_\mH \right)^{-1} \Pi_{H}$, as the pseudoinverse of $\iota$ on $\mH$, in analogy with the Moore-Penrose pseudoinverse of bounded, closed-range linear maps between Hilbert spaces \cite{BeutlerRoot73}. In particular, note that $ T \iota f = f $ for every $ f \in \mH$ and $T g = 0$ for every $g  \in H^\perp$, which shows that $T$ reduces to the Moore-Penrose pseudoinverse $\iota^+$ of $ \iota $ if $\mH$ is a Hilbert space. In that case, the excess generalization error of the target function $f_\mathcal{H} $ in~\eqref{eqFTarget} is due to the component $ (I - \Pi_H ) Z_\tau $ of $ Z_\tau $ in the orthogonal complement of $H$ in $L^2(\mu_X)$. See \ref{appPInv} for additional details on the Moore-Penrose pseudoinverse. 

\paragraph*{Ambient Hilbert space}
Explicit representations of $ T $ depend on the choice of $\mH$, and among the many such possible choices, in KAF we focus on the case where $\mH$ is a finite-dimensional subspace of an ambient Hilbert space $\mathcal{K}$ that $\iota$ compactly embeds into $L^2(\mu_X)$. As $\iota$ is a compact operator between Hilbert spaces, its adjoint $\iota^*:  L^2(\mu_X)  \to \mathcal{K}$ is well-defined and compact. Consequently, the self-adjoint operator $G:= \iota\iota^* : L^2(\mu_X) \to L^2(\mu_X)$ is also compact. The spectral theorem for compact, self-adjoint operators thus guarantees the existence of an orthonormal basis $\{\phi_i\}_{i=1}^\infty$ of $L^2(\mu)$ consisting of eigenfunctions of $G$, with non-negative corresponding eigenvalues $\lambda_i$. 

\begin{rk}
    As we will see in Section~\ref{secRKHS} below, under natural assumptions, $\mathcal K$ has the structure of an RKHS. In that case, the adjoint operator $\iota^*$ becomes an integral operator associated with the reproducing kernel of $ \mathcal K$, and under appropriate continuity assumptions, the orthonormal functions $\psi_i$ correspond to Mercer feature vectors, used, e.g., for unsupervised learning in kernel principal component analysis (KPCA) \cite{ScholkopfEtAl98}. This perspective of expressing integral operators arising in learning problems as adjoints of restriction maps was also adopted by Rosasco et al.\ \cite{RosascoEtAl10} in a study on spectral approximation of integral operators.
\end{rk}

By convention, we order the eigenvalues $\lambda_i$ in decreasing order, so that the sequence $\lambda_1, \lambda_2, \ldots $ only accumulates at zero by compactness of $ G$. Defining

\begin{equation}
    \label{eqPsi}
    \psi_i = \iota^* \phi_i / \lambda_i^{1/2}
\end{equation}
for each $ \lambda_i > 0 $, and choosing  $\ell \in \mathbb{N}$ such that $\lambda_{\ell} >0 $,  we then select as a hypothesis space the $\ell$-dimensional subspace $\mH_\ell \subseteq \mathcal{K}$, where 
\begin{equation}
    \label{eqHypothesis}
    \mH_\ell = \spn \{ \psi_1, \ldots, \psi_\ell \}.
\end{equation}   
It follows from orthonormality of the $ \phi_i $ and their definition in~\eqref{eqPsi} that the $ \psi_i $ form an orthonormal set in $ \mathcal{K} $, i.e., $\langle \psi_i, \psi_j \rangle_\mathcal{K} = \delta_{ij}$. Here, $\langle \cdot, \cdot \rangle_\mathcal{K}$ is the inner product of $\mathcal{K}$, taken conjugate-linear in its first argument. Moreover, the $\psi_i$ are orthonormal eigenfunctions of the operator $ \tilde G := \iota^* \iota $ on $ \mathcal K $, corresponding to the same eigenvalues $\lambda_i$, $\tilde G \psi_i = \lambda_i \psi_i $. In fact, the square roots of the nonzero eigenvalues $\lambda_i$ are the singular values of the compact operator $\iota^*$, and the corresponding $\psi_i$ and $\phi_i$ are left and right singular vectors, respectively, i.e., 
\begin{displaymath}
    \iota^* h = \sum_{i:\lambda_i > 0 } \psi_i \lambda_i^{1/2} \langle \phi_i, h \rangle_{L^2(\mu_X)}, \quad \forall h \in L^2(\mu_X).
\end{displaymath}

With these definitions, it follows that  $\Pi_{H_\ell}$, where $ H_\ell := \iota \mH_\ell $, is the $L^2(\mu_X)$-orthogonal projection with range  $ \spn \{ \phi_1, \dots, \phi_\ell \} $. As for the inverse $(\iota|_{\mH_\ell})^{-1}$, it acts as 
\begin{displaymath}
    (\iota|_{\mH_\ell})^{-1} \phi_i = (\iota|_{\mH_\ell})^{-1} \iota \iota^* \phi_i / \lambda_i =   \psi_i / \lambda_i^{1/2} 
\end{displaymath}
on each eigenfunction $ \phi_i $ corresponding to a nonzero eigenvalue $ \lambda_i $. Consequently, by expanding $Z_\tau \in L^2(\mu_X)$ as $ Z_\tau = \sum_{i=1}^\infty \alpha_i(\tau) \phi_i $,
where 
\begin{equation}
    \alpha_i(\tau) =\langle \phi_i \circ X, U^\tau Y \rangle_{L^2(\mu)},
\label{eq:alpha}
\end{equation}
the target function $f_{\tau,\ell}$ from Definition~\ref{defnTarget} associated with $\mathcal{H}_{\ell}$ is given by
\begin{equation}
    f_{\tau,\ell} = T_\ell Z_\tau =  \sum_{i=1}^\ell \frac{\alpha_i(\tau)}{\lambda_i^{1/2}} \psi_i
\label{eq:theFormula}
\end{equation}
where $T_\ell : L^2(\mu_X) \to \mH_\ell$ is the pseudoinverse operator from~\eqref{eqFTarget}. 

Considering now the image $ K = \iota \mathcal{K} $ of the ambient Hilbert space under $L^2(\mu_X)$ inclusion, one can verify that it can be characterized as the subspace $ K = \left \{ \sum_{i: \lambda_i > 0} c_i \phi_i \in L^2(\mu_X) : \sum_{i:\lambda_i>0} \lvert c_i \rvert^2 / \lambda_i < \infty \right \} $.
The following is then a direct consequence of the definition of the Moore-Penrose pseudoinverse in Definition~\ref{defnPInv}.
\begin{lem}
    \label{lemNystrom}The operator $ \tilde T : D(\tilde T) \to \mathcal{K} $, with dense domain $D( \tilde T) = K \oplus \ker \iota \subseteq L^2(\mu_X)$, defined as $ \tilde T f = \sum_{i:\lambda_i>0} \langle \phi_i, f \rangle_{L^2(\mu_X)} \psi_i / \lambda_i^{1/2} $, is a closed-range operator whose pseudoinverse $\tilde T^+ : \mathcal K \to L^2(\mu_X)$ is equal to $ \iota$. Moreover,  $\tilde T$ is equal to the pseudoinverse of $\iota$, and by~\eqref{eqPInv} we have,
    \begin{displaymath}
        \tilde T f = \iota^* G^+ f, \quad \forall f \in D(G^+).
    \end{displaymath}
\end{lem}

Lemma~\ref{lemNystrom} shows that $ \tilde T $ maps each $L^2(\mu_X) $ equivalence class in its domain to an everywhere-defined function in $ \mathcal{K} $, and whenever $f $ lies in $K$, $ \iota \tilde f = f$. That is, $\tilde T$ is an extension operator, mapping $ f \in K$ to a representative in $\mathcal K$. Note that $ \tilde T $ is necessarily an unbounded operator if $ K $ is infinite-dimensional, and, moreover, if  $G $ is strictly positive-definite (so that all $ \lambda_i $ are strictly positive), then $K$ is a proper, dense subspace of $L^2(\mu_X)$.  In fact, $ \tilde T $ is closely related to the Nystr\"om extension operator employed in applications such as function interpolation and kriging \cite[e.g.,][]{CoifmanLafon06b,RosascoEtAl10}.  Noticing from~\eqref{eq:theFormula} that  $ T_\ell = \tilde T |_{H_\ell} $, we may therefore interpret the target function $f_{\tau,\ell} $ as a spectrally truncated Nystr\"om extension of $ Z_\tau $, which is well defined even if $ Z_\tau $ does not lie in $K$. In fact, it follows from~\eqref{eqFTarget} that $ \iota T_{\ell} $ is equal to the $L^2(\mu_X)$ orthogonal projection $ \Pi_{H_\ell}$. Moreover, we have:  
\begin{lem}
    \label{lemProj}
    As $\ell \to \infty$, $ \iota T_{\ell}$ converges strongly to the orthogonal projection $ \Pi_{\overline K} : L^2(\mu_X) \to L^2(\mu_X)$ onto the $L^2(\mu_X)$-closure of $K$; that is,
    \begin{displaymath}
        \lim_{\ell\to\infty} \iota T_\ell f = \Pi_{\overline K} f, \quad \forall f \in L^2(\mu_X).
    \end{displaymath}
\end{lem}

Lemma~\ref{lemProj} indicates that even if the target functions $ f_{\tau,\ell} = T_{\ell} Z_\tau $ do not have a limit in the ambient space $\mathcal K$, they have an $L^2(\mu_X)$ limit. In particular, if it can be arranged that $K$ is a dense subspace of $L^2(\mu_X)$, $\Pi_{\overline K} = \Id$, and the $f_{\tau,\ell}$ converge in $L^2(\mu_X)$ norm to the regression function $Z_\tau$. Ensuring that $\mathcal K$ is an empirically constructible space with dense image $K = \iota \mathcal K$ in $L^2(\mu_X)$ is a key consideration in KAF, which will occupy us in the ensuing sections.

\subsection{\label{secRKHS}Reproducing kernel Hilbert spaces}

For the remainder of the paper, we will restrict attention to the case that the response variable $Y$ is complex-valued, i.e., $\mathcal{Y} = \mathbb{C}$. In this setting, the ambient Hilbert space $\mathcal{K}$ naturally acquires the structure of an RKHS \cite{Paulsen16,FerreiraMenegatto13}, as we describe below.

\begin{defn}[RKHS] For each point $x \in \mathcal X$, let $L_x : \mathcal{K} \to \mathbb{C}$ be the evaluation functional on the ambient Hilbert space, defined by  $L_x f = f(x)$. The space $\mathcal{K}$ is said to be an RKHS if $L_x$ is bounded, and therefore continuous, at every $x \in \mX$.  
\end{defn}
    
It is a known fact that no unbounded linear functional on a Banach space can be constructed without the axiom choice. Therefore, all explicitly constructible Hilbert spaces of complex-valued functions are necessarily RKHSs. Consequently, all explicitly representable target functions $f_{\tau,\mH}$ from~\eqref{eqFTarget} necessarily lie in an RKHS. Note that by boundedness of $L_x $ at every $ x \in \mX$, convergence of two functions in RKHS norm implies pointwise convergence on $\mX$.

\paragraph*{Basic properties of RKHSs} It follows from the Riesz representation theorem that for every $x \in \mX$  there exists some function $k_x \in \mathcal{K}$ such that
\begin{displaymath}
f(x) = L_x f = \langle k_x, f \rangle_\mathcal{K}, \quad \forall f \in \mathcal{K}. 
\end{displaymath}
The above is known as the reproducing property of $\mathcal{K}$. The \emph{reproducing kernel} $k : \mX \times \mX \to \mathbb{C}$ of $\mathcal{K}$ is then defined as the bivariate function 
\begin{displaymath}
k(x_1, x_2)  = \langle k_{x_1}, k_{x_2} \rangle_\mathcal{K}.
\end{displaymath}
It follows from the defining properties of inner products that $k$ is (i) conjugate-symmetric, i.e., $ k( x_1, x_2 ) = k(x_2,x_1)^* $ for all $x_1, x_2 \in \mX$; and (ii) positive-definite, i.e., for all $ x_1,\ldots, x_m \in \mX $ and $a_1,\ldots,a_m \in \mathbb{C}$, 
\begin{equation}
    \label{eqKPosDef}
\sum_{i=1}^m \sum_{j=1}^m a_i^* a_j  k(x_i, x_j) \ge 0.
\end{equation}
Conversely, the Moore-Aronszajn theorem \cite{Aronszajn50} states that for any conjugate-symmetric, positive-definite kernel function $k : \mX \times \mX \to \mathbb{C}$, there exists a unique RKHS on $\mX$ for which $k$ is the reproducing kernel. Thus, there is a one-to-one correspondence between kernels and RKHSs. 

Let $\rho : \Sigma_\mX \to [ 0, \infty ] $ be any measure such that there exists a compact embedding $\iota_\rho$ of the RKHS $\mathcal{K}$ into $L^2(\rho)$. The practical utility of RKHSs manifests in the adjoint $\iota^*_\rho : L^2(\mX) \to \mathcal{K}$ being representable in terms of the kernel as
\begin{displaymath}
    \iota^*_\rho f(x) = \langle k_x, \iota_\rho^* f \rangle_\mathcal{K} = \langle \iota_\rho k_x, f \rangle_{L^2(\rho)} = \int_\mX k^*( x, \cdot ) g \, d\rho ,
\end{displaymath}
where $f$ is any element of  $L^2(\rho)$. Thus, the adjoint of the embedding of the RKHS $\mathcal{K} $ into $L^2(\rho)$ is a compact integral operator on the latter space. Similarly, $ G_\rho := \iota_\rho \iota^*_\rho $ is a positive-semidefinite, self-adjoint, compact integral operator on $L^2(\rho)$. 

\begin{defn}[particular classes of kernels] \label{defnKernels} We will say that a positive-definite kernel $k : \mathcal X \to \mathcal X \to \mathbb C $ is:
    \begin{itemize}[wide]
        \item Strictly positive-definite if the inequality in~\eqref{eqKPosDef} is strict whenever the $x_i $ are all distinct and at least one of the $a_i $ is nonzero.
    
        \item $L^2(\rho)$-strictly-positive   if $G_\rho$ is a strictly-positive operator. In that case, $ K_\rho := \iota_\rho \mathcal{K}$ is a dense subspace of $L^2(\rho)$. 
    
        \item $L^2(\rho)$-Markov if $G_\rho $ is a Markov operator, i.e., (i) $\int_\mX G_\rho f \, d\rho = \int_\mX f \, d\rho$ for all $ f \in L^2(\rho)$; and (ii) $G_\rho f = f$ if $f$ is constant. As a result, the leading largest eigenvalue of $G_\rho $ is equal to 1, and the corresponding eigenspace contains constant functions. 

        \item An $L^2(\rho)$-Markov ergodic if it is $L^2(\rho)$-Markov and the eigenvalue 1 of $G_\rho$ is simple.
    \end{itemize}
\end{defn}

In the case $\rho = \mu_X$, we will abbreviate $ \iota_{\mu_X} = \iota $ and $G_{\mu_X} = G$ as in Section~\ref{secHypothesis}. The evaluation of the target function from~\eqref{eqFTarget} at a point $x \in \mX$ is then expressible as
\begin{equation}
    f_{\tau,\ell}(x) = \sum_{i=1}^\ell \frac{\alpha_i(\tau)}{\lambda_i} \langle \iota k_x, \phi_i \rangle_{L^2(\mu_X)} = \sum_{i=1}^\ell \frac{\alpha_i(\tau)}{\lambda_i^{1/2}} \psi_i(x).
\label{eq:firstKPCR}
\end{equation}

\paragraph*{Topological framework and Mercer kernels}
Henceforth, we will assume that $\mX$ has the structure of a metric space, equipped with its Borel $\sigma$-algebra $\Sigma_\mX$, and $\mu_X$ is a Borel probability measure with compact support $ \mX_\mu \subseteq \mX$. Given any subset $S$ of $\mX$, we will use the notation $\mathcal{K}(S)$ to represent the RKHS on $S$ with reproducing kernel $k|_{S\times S}$. Note that $\mathcal{K}(S)$ embeds naturally and isometrically into $\mathcal{K}$, so we may view it as a subspace of the latter space. We also let $C(S)$ be the space of complex-valued continuous functions on $S$, and $C_b(S)$ the Banach space of bounded functions in $C(S)$, equipped with the uniform norm. Note that $C(S) = C_b(S) $ if $S $ is compact.

In this setting, continuous kernel functions on $\mX$, also known as Mercer kernels, have the property that their associated RKHS is a subset of $C(\mX)$ \cite{FerreiraMenegatto13}. Moreover, for any compact set $ S \subseteq \mX $, the embedding $\mathcal{K}(S) \hookrightarrow C(S)$ is bounded. If, in addition, $S$ is the support $\mX_\rho$ of a finite Borel measure $ \rho $ on $ \mX$,  $C(\mX_\rho) $ embeds into $L^2(\rho)$ via a bounded linear map, and thus  $\iota_\rho : \mathcal{K}(\mX_\rho) \to L^2(\rho)$ is a bounded, injective operator. It also follows by continuity of $ k $ and compactness of $ \mX_\rho $ that $G_\rho= \iota_\rho \iota^*_\rho$ is a trace-class (and therefore compact) operator, with trace norm equal to $ \tr G_\rho = \int_\mX k(x,x) \, d\rho(x)$ \cite{Brislawn91}. In particular, the compactness of $G_\rho$ is equivalent to $ \iota_\rho $ being compact. Mercer's theorem \cite[][Section~11.4]{Paulsen16} also states that for any $x,x' \in X_\rho$ the kernel $k(x,x') $ can be expressed through the series expansion, 
\begin{equation}
    \label{eqMercer}
    k(x,x') = \sum_{i: \lambda_i > 0} \psi^*_i(x) \psi_i(x'),
\end{equation}
where the $ \psi_i $ are orthonormal functions in $ \mathcal{K} $ associated with eigenvalue $ \lambda_i $ of $ G_\rho $, defined analogously to~\eqref{eqPsi}, and convergence of the sum over $ i$ is uniform on $ \mathcal{X}_\rho \times \mathcal{X}_\rho$. This result then implies that the restrictions of the $ \psi_i $ on $X_\rho$ form an orthonormal basis of $\mathcal{K}(X_\rho)$ (as opposed to merely an orthonormal set). It can also be shown that every strictly positive-definite Mercer kernel is $L^2(\rho)$-strictly positive for any compactly supported, finite Borel measure $ \rho $. See \cite{SriperumbudurEtAl11} for a detailed study on the relationships between the RKHSs associated with different kernel classes (including those in Definition~\ref{defnKernels}) and spaces of functions and measures, such as spaces of continuous functions and $L^p$ spaces. 

By virtue of the above properties, Mercer kernels provide a convenient practical means of generating hypothesis spaces that are compactly embedable into $L^2(\mu_X)$, as required for the hypothesis spaces in Section~\ref{secHypothesis}. Note that the target function in \eqref{eq:firstKPCR} associated with a Mercer kernel is an RKHS (and thus continuous) function defined on the whole of $\mX$, but its behavior outside of the support $ \mX_\mu$ makes no contribution to the excess generalization error from~\eqref{eqExcessGen} determined through the $L^2(\mu_X)$ norm. 

\begin{rk}
    The Mercer expansion in~\eqref{eqMercer} allows one to evaluate inner products between distinguished elements of the RKHS, namely the kernel sections $k_x$ simply by evaluation of a known kernel function, $\langle k_x, k_{x'} \rangle_{\mathcal K} = k(x,x')$ (i.e., the left-hand side of~\eqref{eqMercer}), without having to compute a potentially infinite set of basis vectors for the space (i.e., the eigenfunctions $\psi_i$ in the right-hand side). This well known ``kernel trick'' is employed in a variety of learning techniques, including kernel KRR and SVMs \cite{SteinwartChristmann08}. In contrast, in KAF/KPCR methods incur a potentially significant computational cost associated with computing (training phase) and evaluating (prediction phase) a set of orthogonal basis functions $ \{ \psi_1, \ldots, \psi_{\ell} \} $ with $ \ell \ll n$, with the benefit of controlling the regularity of the target functions through the spectral truncation parameter $\ell$. As we will see in Section~\ref{secSampleError} below, this is an effective means of controlling the sample error, allowing the method to operate stably in environments with small training datasets.
\end{rk}

\subsection{\label{secTarget}Data-driven target function}

We are now ready to construct the empirical target function employed in KAF. In this construction we consider a standard supervised learning scenario, where we have access to a training dataset consisting of pairs $(x_1,y_{1}), (x_2, y_{2}), \ldots, (x_n, y_{n})$, where $ x_j = X(\omega_j)$ and $ y_{j} =  Y(\omega_j)$ are the values of the covariate and response variables, respectively, on an (unknown) collection of points $\omega_1, \ldots, \omega_n$ on the sate space $ \Omega $, sampled along a single dynamical trajectory 
\begin{equation}
    \label{eqTraj}
    \omega_j = \Phi^{t_j}(\omega_1), \quad t_j = (j-1) \, \Delta t, 
\end{equation}
at a fixed sampling interval $\Delta t > 0$. Alternatively, the $(x_j, y_{j} )$ may be generated by an ensemble of (shorter) trajectories on $\Omega$, so long as the timespan of each of these trajectories is not smaller than the desired lead time $\tau$.

\paragraph{Sampling measures} Associated with every dataset from~\eqref{eqTraj} is an empirical probability measure $\mu_n : \mathcal F \to [0,1]$, defined as $ \mu_n = \sum_{j=1}^n \delta_{\omega_j} / n $,
where $\delta_{\omega_j}$ is the Dirac $\delta$-measure supported on $\{ \omega_j \} \subset \Omega $. Similarly, the empirical probability measure $\mu_{X, n} : \Sigma_\mathcal{X} \to [0,1]$ is defined as $\mu_{X,n} =  \sum_{j=1}^n \delta_{x_j} / n $. Intuitively, we view $\mu_n$ and $\mu_{X,n}$ as empirical approximations to  $\mu$ and $ \mu_X$, respectively; a connection which will be made precise in Section~\ref{secError}. 

Next, as empirical analogs of $L^2(\mu)$ and $L^2(\mu_X)$, we employ the Hilbert spaces $L^2(\mu_n)$ and $ L^2(\mu_{X,n})$, consisting of equivalence classes of complex-valued, measurable functions on $\Omega $ and $X$ having common values at the sampled points $\omega_j$ and $x_j$, respectively. As Hilbert spaces, $L^2(\mu_n)$ and $L^2(\mu_{X,n})$ have dimension at most $n$ (with equality if all $\omega_j$ and $ x_j $ are distinct, respectively), and can be homomorphically embedded into $\mathbb{C}^n$, equipped with the normalized dot product $ \bm f \cdot \bm g / n$. That is, for every measurable function $ f : \Omega \to \mathbb{C}$, the corresponding $L^2(\mu_n)$ equivalence class can be represented by a column vector $ \bm f \in L^2(\mu_n)$ with $\bm f = [ f(\omega_1), \ldots, f(\omega_n)]^T $, storing in its components the values of $f $ on $\omega_j$. Elements of $L^2(\mu_{X,n})$ are represented by $\mathbb{C}^n$ vectors in a similar manner, while operators on $L^2(\mu_n)$ and $L^2(\mu_{X,n})$ are represented by $n \times n$ complex matrices. 

As in the case of the $L^2(\mu)$ and $L^2(\mu_X)$ spaces, there is an isometric embedding $ \Xi_n : L^2(\mu_{X,n})\to L^2(\mu_n)$, given by composition by the covariate map, $ \Xi_n f = f \circ X $, whose image we denote by $ L^2_X(\mu_n) = \Xi_n L^2(\mu_{X,n}) $. We also let $ \Pi_{X,n} : L^2(\mu_n) \to L^2(\mu_n)$ be the orthogonal projection mapping into $L^2_X(\mu_n)$. Note that in typical applications involving distinct training data, $ \Pi_{X,n}$ is the identity map and $ \Xi_n$ is unitary, even if $ X $ is non-injective on sets of positive $\mu$-measure (in which case, $\Pi_X$ is not the identity). This disparity between $ \Pi_X $ and $ \Pi_{X,n}$ highlights the risk of overfitting commonly faced by data-driven techniques, which KAF addresses by employing hypothesis spaces of significantly lower dimension than the number of training samples.  

\paragraph{Shift operators} In order to parallel the construction of the target function $ f_{\tau,\ell}$ from Section~\ref{secHypothesis}, we would now like to define a Koopman operator on $L^2(\mu_n)$. However, an obstruction to this is that, unlike the $L^2(\mu)$ setting associated with the invariant measure, the composition operator with respect to the dynamical flow does not lift to an operator on equivalence classes of functions on the $L^2(\mu_n)$ spaces associated with the sampling measures. This is because the flow $\Phi^\tau : \Omega \to \Omega$ on state space does not preserve null sets with respect to $ \mu_n$, meaning that if $ f, f' : \Omega \to \mathbb C$ are measurable functions lying in the same $L^2(\mu_n)$ equivalence class, their images $ U^\tau f = f \circ \Phi^t $ and $U^\tau f' = f' \circ \Phi^t$ may lie in different $L^2(\mu_n)$ equivalence classes. Nevertheless, for any $ q \in \mathbb N_0$, an analogous notion  to the Koopman operator $ U^{q \, \Delta t}$ on $L^2(\mu)$ is provided by the shift operator $U^{q}_n : L^2(\mu_n) \to L^2(\mu_n)$ \cite{BerryEtAl15}, defined as 
\begin{displaymath}
    U^{q}_n f( \omega_j ) = 
    \begin{cases} 
        f(\omega_{j+q}), & j + q \leq n, \\
        0, & \text{otherwise}.
    \end{cases}
\end{displaymath}
Hereafter, we will refer to the $ \mathbb{C}^n$ vector 
\begin{equation}
    \label{eqAnalogVec}
    \bm y_\tau =  [y_{1+q}, \dots, y_{n},0, \dots, 0]^T,
\end{equation}
representing the response $ U^{q}_n \iota_n Y \in L^2(\mu_n)$ for $ \tau = q \, \Delta t$, as the \emph{analog vector}. 

\paragraph{Empirical error minimization} With these definitions, and assuming throughout that $ \tau = q \, \Delta t$ for some $ q \in \mathbb N_0$, the empirical generalization error $\mE_{\tau,n} : \BL^2(\mu_{X,n}) \to \BR$ is given by (cf.~\eqref{eqGen})
\begin{displaymath}
\mE_{\tau,n}(f) : =\lVert \iota_n f \circ X - U^\tau Y \lVert^2_{L^2(\mu_n)},
\end{displaymath}
where $ \iota_n : \BL^2(\mu_{n}) \to L^2(\mu_{n})$ maps each function in $\BL^2(\mu_n)$ to its $L^2(\mu_n)$ equivalence class. This functional is minimized by a unique element $Z_{\tau,n} \in L^2(\mu_{X,n}) $ analogous to the regression function $Z_\tau$ from Sections~\ref{secBackground}--\ref{secRKHS}. Moreover, we may split the empirical generalization error as 
\begin{displaymath}
    \mE_{\tau,n}(f) = \mA_{\tau,n}(f) + \sigma_{\tau,n}, 
\end{displaymath}
with (cf.~\eqref{eqExcessGen}) 
\begin{displaymath}
    \mA_{\tau,n}(f) = \lVert \iota_n f - Z_{\tau,n} \lVert^2_{L^2(\mu_{X,n})}, \;\; \sigma_{\tau,n} = \lVert Z_{\tau,n} \circ X - U^\tau Y \rVert_{L^2(\mu_n)}^2.
\end{displaymath}

\paragraph{Empirical hypothesis space} To construct an empirical target function, we proceed again analogously to the infinite-dimensional case in Sections~\ref{secBackground}--\ref{secRKHS}. That is, we seek the minimizer of the empirical excess generalization error $\mA_{\tau,n}(f) $ for $f$ lying in an $\ell$-dimensional empirical hypothesis space $\mH_{\ell,n}$, which is chosen as a subspace of an ambient RKHS $\mathcal{K}_n \subset C(\mX) $ associated with an empirical Mercer kernel $k_n : \mX \times \mX \to \mathbb{C}$. Note that we allow the reproducing kernel $k_n$ to depend on $n $ in order to be able to take advantage of the variety of normalized kernel algorithms in the literature \cite{CoifmanLafon06,VonLuxburgEtAl08,CoifmanHirn13,BerryHarlim16,BerrySauer16b}. Given any $x \in \mX$, we shall refer to the $\mathbb{C}^n$ vector  
\begin{displaymath}
    \bm k(x) = [k_n(x, x_1), \dots, k_n(x, x_n)]^T,
\end{displaymath}
representing the $L^2(\mu_{X,n})$ equivalence class $\iota_n k_n(x,\cdot)$ of the kernel section $k_n(x,\cdot) \in \mathcal{K}_n$ as the \emph{kernel vector}.

Next, because $ \mathcal{K}_n \subseteq \mathbb{L}^2(\mu_{X,n})$, we can consider $\iota_n : \mathcal{K}_n \to L^2(\mu_{X,n})$ as a (finite-rank, and thus compact) operator between Hilbert spaces, inducing the self-adjoint integral operator $ G_n := G_{\mu_{X,n}} = \iota_n \iota_n^*$ on $L^2(\mu_{X,n})$. The leading $\ell$ orthonormal eigenvectors $ \phi_{1,n}, \ldots,  \phi_{n,\ell} $ of $G_n$, corresponding to positive eigenvalues $\lambda_{1,n} \geq \cdots \geq \lambda_{\ell,n} $, respectively, induce the $\ell$-dimensional empirical hypothesis space $ \mH_{\ell,n} \subseteq \mathcal{K}_n$ given by (cf. \eqref{eqHypothesis}) 
\begin{displaymath}
    \mH_{\ell,n} = \spn \{ \psi_{1,n}, \ldots, \psi_{\ell,n} \},
\end{displaymath}
where
\begin{equation}
    \label{eqPsiN}
    \psi_{i,n} = \iota_n^* \phi_{i,n} / \lambda_{i,n}^{1/2}
\end{equation}
are orthonormal functions in $ \mathcal{K}_n$.  We then compute the minimizer$f_{\tau,\ell,n} \in \mH_{\ell,n}$ of $\mA_{\tau,n}(f)$ over this hypothesis space, obtaining, in direct analogy to~\eqref{eq:firstKPCR}, 
\begin{align}
    \nonumber
    f_{\tau,\ell, n}(x) = T_{\ell,n} Z_{\tau,n} &=  \sum_{i=1}^\ell \frac{\alpha_{i,n}(\tau)}{ \lambda_{i,n}} \left\langle \iota_n k_{n,x},  \phi_{i,n}\right\rangle_{L^2(\mu_{X,n})} \\
    &= \sum_{i=1}^{\ell} \frac{\alpha_{i,n}(\tau)}{\lambda_{i,n}^{1/2}}\psi_{i,n}(x),
    \label{eqKAF}
\end{align}
where $ \alpha_{i,n}(\tau) = \langle \phi_{i,n}, U^{q}_n Y \rangle_{L^2(\mu_{X,n})}$, and $T_{\ell,n} : L^2(\mu_{X,n}) \to \mH_{\ell,n} $ is the pseudoinverse operator associated with $\mH_{\ell,n}$. This minimizer constitutes the empirical target function utilized by KAF. 

\begin{defn}[empirical target function]\label{defnEmpTarget}The empirical target function $f_{\tau,\ell,n}$ at lead time $\tau $ associated with the response $Y$ and $\ell$-dimensional hypothesis space $\mathcal H_{\ell,n}$ is the minimizer of the empirical excess generalization error functional $\mathcal A_{\tau,n}$ over $\mathcal H_{\ell,n}$, given by~\eqref{eqKAF}.
\end{defn}

The expression in~\eqref{eqKAF} can be written more compactly in matrix form using the column vector representations $ \bm \phi_i \in \mathbb{C}^n $ of the $ \phi_{i,n} $, given by eigenvectors of the $n \times n $ kernel matrix $\bm G = [ k^*_n(x_i,x_j) ] /n$ representing $G_n $, and chosen such that $ \bm \phi_i \cdot \bm \phi_j / n = \delta_{ij}$. Note, in particular, that the expansion coefficients $ \alpha_{i,n}(\tau) $ are simply equal to the dot products $ \alpha_{i,n}(\tau) = \bm \phi_i \cdot \bm y_\tau / n $ with the analog vector. Treating the remaining terms in~\eqref{eqKAF} in a similar manner, we arrive at the expression     
\begin{equation}
    \label{eqAMat}
    f_{\tau,\ell,n}(x) = \bm k (x)^{*} \bm A_\ell  \bm y_\tau, \quad \bm A_\ell = \bm \Phi_\ell \bm \Lambda^{-1}_\ell \bm \Phi^{*}_\ell / n^2,
\end{equation}
where $ \bm k( x ) $ is the kernel vector, $\bm \Phi_\ell$ is the $n \times \ell$ matrix whose columns consist of the eigenvectors  $\bm \phi_i \in \mathbb{C}^n $, $\bm \Lambda_\ell $ is the $\ell \times \ell$ diagonal matrix whose diagonal entries consist of $\lambda_{i,n}$, and $^*$ denotes complex-conjugate transpose. This formula expresses the KAF target function as a sesquilinear form $( \bm k(x), \bm y_\tau) \mapsto \bm k (x) ^{*}   \bm A_\ell \bm y_\tau $, mapping pairs of kernel and analog vectors to $\mathbb{C}$-valued forecasts.  Letting $\bm V_\ell = \bm \Lambda^{-1/2}_\ell \bm \Phi^*_\ell/ n$, where $  \bm V_\ell{}^* \bm V_\ell = \bm A_\ell$, the empirical target function is reexpressed as
\begin{displaymath}
f_{\tau,\ell,n} (x) =  \bm V_\ell \bk(x) \cdot \bm V_\ell \bm y_\tau .
\end{displaymath}
This particular form emphasizes that the forecast is the result of taking the inner product of suitably projected kernel vector and equivalently projected analog vectors. 

\begin{rk}
    In KPCA \cite{ScholkopfEtAl98}, as well as related manifold learning techniques \cite{BelkinNiyogi03,HeinEtAl05,CoifmanLafon06,Singer06,BerryHarlim16,BerrySauer16b}, eigenvectors of kernel matrices such as $ \bm \phi_i $ above are employed for unsupervised feature extraction. In particular, it is common to use the $ \bm \phi_i $ as coordinate vectors of  dimension reduction maps, $ x_j \mapsto ( \lambda_1^{-1/2} \phi_{1,n}( x_j ), \ldots, \lambda_\ell^{-1/2} \phi_{\ell,n}(x_i) ) \in \mathbb C^\ell $, mapping potentially high-dimensional covariate data into low-dimensional Euclidean spaces, where geometrical data relationships are revealed. In contrast, KAF/KPCR is a supervised learning technique where the goal is to perform out-of-sample prediction of a random variable (the response observable $U^\tau Y$). A common aspect of the two methods is that they both rely heavily on eigendecompositions of kernel integral operators, but the end goals are fundamentally different. Note, in particular, that in KPCR one seeks to use as many eigenvectors $ \bm \phi_i $ that can be computed from the available training data with tolerable sample error, and the number $\ell$ of such eigenvectors can be far higher than the dimension of covariate space $\mathcal X$. For instance, in the  L63 examples in Figure~\ref{figL63Pred} and Section~\ref{secL63} we use $\ell = O(10^3) $, which clearly does not serve as a ``dimension reduction'' map for the 3- or 1-dimensional covariate spaces employed there.   
\end{rk}

\section{\label{secError}Error analysis and convergence}

The previous section has shown how to calculate both an empirical target function $f_{\tau,\ell,n}$ (Definition~\ref{defnEmpTarget}) and an ideal target function $f_{\tau,\ell}$ (Definition~\ref{defnTarget}), corresponding to two different hypothesis spaces, $\mH_{\ell,n}$ and $\mH_{\ell}$, as well as two different error functionals, $\mE_{\tau,n}$ and $\mE_{\tau}$, respectively. This section addresses the connection between the two functions, with the ultimate goal being that of bounding the error $\mE_\tau(f_{\tau,\ell,n})$ of the empirical target function as much as possible. Among other reasons, the availability of such bounds is useful for assessing the risk of overfitting the training data; that is, the possibility that $\mE_\tau(f_{\tau,\ell,n}) \gg \mE_{\tau,n}(f_{\tau,\ell,n})$ for the chosen empirical hypothesis space. Note, in particular, that for a variety of kernels (e.g., strictly positive-definite kernels) it is possible to make $\mE_{\tau,n}(f_{\tau,\ell,n})$ at fixed $n$ arbitrarily small by increasing $\ell$, but this reduction of empirical error eventually leads to an increase of the ``true'' error $\mE_\tau(f_{\tau,\ell,n})$ with respect to the invariant measure of the dynamics. See Section~\ref{secCircle} for an illustration of this phenomenon. 

The analysis of the error $\mE_\tau(f_{\tau,\ell,n})$ is typically organized into analysis of the error $\mE_\tau(f_{\tau,\ell})$ of the ideal target function (i.e, the excess generalization error), and the difference in error $\mE_\tau(f_{\tau,\ell,n}) - \mE_\tau(f_{\tau,\ell})$, denoted by $\mathcal{D}_{\tau,\ell,n}$, between the empirical and ideal target functions, referred to as the sample error \cite{cucker2007learning}.   In other words, error analysis uses the following decomposition:
\begin{align*}
    \mE_\tau(f_{\tau,\ell,n})  &= \mE_\tau(f_{\tau,\ell}) + \mathcal{D}_{\tau,\ell\,n} \\
    &= \mE_\tau(f_{\tau,\ell}) + \left(\mE_\tau(f_{\tau,\ell,n}) - \mE_\tau(f_{\tau,\ell})\right).
\end{align*}
This section examines in detail these contributions, and establishes sufficient conditions for convergence of the KAF target function to the conditional expectation.  

\subsection{\label{secGenError}KAF generalization error}

The excess generalization error $\mA_\tau(f_{\tau,\ell})$ from~\eqref{eqExcessErrFl} of the KAF target function $f_{\tau,\ell}$ in~\eqref{eq:firstKPCR} is given by
\begin{displaymath}
\mA_\tau(f_{\tau,\ell}) =  \lVert \Pi_{H_\ell^{\perp}} Z_\tau \rVert^2_{L^2(\mu_X)} = \sum_{i=\ell+1}^\infty \lvert \alpha_i(\tau) \rvert^2,
\end{displaymath}
where $\Pi_{H_\ell^{\perp}} : L^2(\mu_X) \to L^2(\mu_X)$ is the orthogonal projection mapping into the orthogonal complement  $H_\ell^\perp$ of the hypothesis space $H_\ell$ in $L^2(\mu_X)$. It follows from the above that $\mA_\tau(f_{\tau,\ell}) $ vanishes as $ \ell \to \infty $ for any $Z_\tau \in L^2(\mu_X) $, and thus for any response variable $U^\tau Y \in L^2(\mu) $, iff the sequence of projections $ \Pi_{H_\ell^\perp} $ converges pointwise to 0 as $ \ell \to \infty $ (i.e., $ \Pi_{H_\ell^\perp} g \to 0 $ for any $ g \in L^2(\mu_X)$). By Lemma~\ref{lemProj}, this happens in turn iff $K$ is a dense subspace of $L^2(\mu_X)$, i.e., iff $ G $ is a strictly positive operator. Since $ \mA_\tau(f) = 0 $ iff $ \iota f = Z_\tau $, we obtain the following basic consistency result expressed in terms of a positivity condition on the kernel $k$. 

\begin{thm}
    \label{thmConv} Let $k: \mX \times \mX \to \mathbb{C} $ be an $L^2(\mu_X)$-strictly-positive kernel with corresponding RKHS $\mathcal{K}$. Then, for any response variable $U^\tau Y \in L^2(\mu)$ and lead time $ \tau \geq 0$, as $\ell\to\infty$, the target functions $f_{\tau,\ell}$ from~\eqref{eq:firstKPCR} converge to the conditional expectation $\mathbb{E}[U^\tau Y \mid X ] = X \circ Z_\tau$, in the sense that $\lim_{\ell \to \infty} \lVert \iota f_{\tau,\ell} - Z_\tau \rVert^2_{L^2(\mu_X)} = 0 $. 
\end{thm}

Convergence with respect to the (stronger) RKHS norm of $\mathcal{K}$, as well as more precise estimates of the $L^2(\mu_X)$ error, can be obtained under the additional assumption that the regression function $Z_\tau$ lies in the subspace $K \subset L^2(\mu_X) $. In that case, $ Z_\tau $ has a representative $ f_\tau \in \mathcal{K} $, given by Nystr\"om extension as 
\begin{displaymath}
    f_\tau = \tilde T Z_\tau = \sum_{i=1}^\infty \alpha_i(\tau) \psi_i / \lambda_i^{1/2},   
\end{displaymath}
where the infinite sum in the right-hand side converges in $\mathcal{K}$ norm. That is, $ f_\tau $ is given by the $\mathcal{K} $-norm limit of the partial sums $ \sum_{i=1}^\ell \alpha_i  \psi_i / \lambda_i^{1/2} $. The latter are precisely equal to the target functions $ f_{\tau,\ell} $ from~\eqref{eq:firstKPCR}, and therefore we conclude that $ \lim_{\ell\to\infty}\lVert f_{\tau,\ell} - f_\tau \rVert_{\mathcal{K}} = 0$.  

To obtain an estimate of $\mathcal{A}(f_{\tau,\ell})$, observe that $K$ coincides with the range of $G^{1/2}$, the square root of $G$. It then follows that for $ Z_\tau \in K $, there exists $ W_\tau \in L^2(\mu_X) $ such that $ Z_\tau = G^{1/2} W_\tau$, which allows the excess generalization error to be rewritten as
\begin{displaymath}
\mA_\tau(f_{\tau,\ell}) = \sum_{i = \ell+1}^\infty \lambda_i \lvert \langle \phi_i,  W_\tau \rangle_{L^2(\mu_X)} \rvert^2.
\end{displaymath}
The Cauchy-Schwarz inequality then yields
\begin{displaymath}
\mA_\tau(f_{\tau,\ell})  \le  \left(\sum_{i=\ell+1}^\infty \lambda_i \right) \lVert W_\tau \rVert^2_{L^2(\mu_X)},
\end{displaymath}
where $ \sum_{i=\ell+1}^\infty \lambda_i < \tr G$ is finite. Thus, in this case we can bound the decay of the excess generalization error by the decay of the tail sum of the eigenvalues of $G$. 

The study of decay rates of the eigenvalues of an integral operator \cite{Konig86} is an active field of research. In the setting of Mercer kernels and compactly supported probability measures studied here, it can be shown that $\lambda_i = o(i^{-1})$ for large-enough $i$ \citep[][Theorem~2.4]{FerreiraMenegatto09}, which is consistent with the fact that $G $ is trace-class. Estimates of the rate of decay of the tail sum are possible under additional regularity conditions on the kernel, including, for example, specialized notions of Lipschitz continuity. In such cases, it is possible to express the decay rate of $\mA_\tau(f_{\tau,\ell})$ as being algebraic, i.e., $\mA_\tau(f_{\tau,\ell})\le C \ell^{-\gamma}$, for some positive constants $C$ and $\gamma$ \cite{ferreira2013eigenvalue}.

\subsection{\label{secSampleError}KAF sample error}

In this section, we will establish that, under natural assumptions on the dynamical system and the reproducing kernels, the difference in error $\mathcal{D}_{\tau,\ell,n}$ between the empirical and ideal target functions vanishes in the limit of large data, $n \to \infty$. We will do so by establishing a stronger result, namely that $f_{\tau,\ell,n}$ converges uniformly to $ f_{\tau,\ell} $ in an appropriate compact set containing the supports of $ \mu_X $ and the sampling measures $\mu_{X,n}$.

\paragraph*{Basic assumptions for convergence} Our first assumption is that (i) $\Omega$ has the structure of a metric space, equipped with its Borel $\sigma$-algebra $ \mathcal F$; (ii) $ \mu $ is a Borel probability measure with compact support $ \Omega_\mu \subseteq \Omega$; and (iii) all of $ \Phi^t $, $X$, and $Y$ are continuous. Note that, by continuity of $ \Phi^t$, the Koopman operator $U^t $ maps continuous functions to continuous functions for all $ t \geq 0 $, preserving the supremum norm of bounded continuous functions in $C_b(\Omega)$. See the commutative diagram in Figure~\ref{figCommut} for an illustration of the relationships between the Koopman operator on $C(\Omega)$ and $L^2(\mu_X)$.    

Our second assumption pertains to the convergence of the empirical measures $ \mu_{n}$ underlying the data to the invariant measure. Specifically, we assume that, for the starting state $\omega_1 \in \Omega$, the measures $\mu_{n}$ converge to $ \mu$ weakly; that is, for every bounded, continuous function $ g: \Omega \to \mathbb{C}$, $\lim_{n\to\infty} \int_\Omega g \, d\mu_{n} = \int_\Omega g \, d\mu$. The weak convergence of $ \mu_n $ to $ \mu $, in conjunction with the continuity of $X$, implies in turn that $\mu_{X,n}$ converges weakly to $ \mu_X$, i.e., $\lim_{n\to\infty} \int_\mX f \, d \mu_{X,n} = \int_X f \, d\mu_X $, for all $f \in C_b(\mX)$. 

Our third assumption relates to the existence of a compact set in which both the the covariate data $x_i $ and the support $ \mu_{X}$ lie. Specifically, for the starting state $\omega_1 \in \Omega$ underlying the covariate training data, we assume that there exists a compact set $\mathcal{U} \subseteq \mX$ containing $\mathcal{X}_\mu = \supp \mu_X $, as well as $ \supp \mu_{X,n} = \{ x_1, \ldots, x_n \} $ for every $ n \in \mathbb{N}$. This condition is automatically satisfied if the state space $ \Omega $ is already a compact space (e.g., ergodic dynamics on a torus), and is also satisfied by many systems with appropriate dissipative dynamics. Examples of such systems include ordinary differential equation models on $\Omega = \mathbb{R}^d$ with quadratic nonlinearities, such as the L63 system \cite{LawEtAl14} studied in Section~\ref{secL63} below, as well as partial differential equation models possessing inertial manifolds \cite{ConstantinEtAl89}. For our purposes, the existence of the compact set $\mathcal{U}$ allows the (generally distinct) ideal and empirical RKHSs, $\mathcal{K}(\mathcal{U})$ and $\mathcal{K}_{n}(\mathcal{U})$, respectively, to be viewed as subspaces of the Banach space $C(\mathcal{U})$. In the latter space, the relevant notion of convergence is convergence with respect to the uniform norm. 

Next, we make an assumption on the convergence of the empirical reproducing kernels $k_n$ of $ \mathcal{K}_n  $ to the reproducing kernel $k$ of $ \mathcal{K}$. Specifically, we assume that, as $n \to \infty$, $k_n $ converges to $ k $ uniformly on $ \mathcal{U} \times \mathcal{U} $ (i.e., with respect to $C(\mathcal{U}\times\mathcal{U})$ norm). This assumption is trivially satisfied if one works with data-independent kernels, $k_n = k $, and also holds for many classes of normalized kernels, including  \cite{HeinEtAl05,CoifmanLafon06,VonLuxburgEtAl08,CoifmanHirn13,BerryHarlim16,BerrySauer16b}.

Finally, we assume that the response variable $U^\tau Y $ is bounded on $ X^{-1}(\mathcal{U})$, i.e., $C_Y = \sup_{\omega \in X^{-1}(\mathcal{U})} \lvert U^\tau Y(\omega) \rvert < \infty $.

\paragraph{Physical measures} We define the basin of $\mu$ as the maximal set $M_\mu \subseteq \Omega $ for which the sampling measures $ \mu_n $, starting from any $ \omega_1 \in M_\mu $, converge weakly to $ \mu $. If the dynamics is ergodic (i.e., every invariant set $S \in \mathcal F$ under $\Phi^t$ for all $t \in \mathbb{R}$ has either $\mu(S) = 0$ or $ \mu(S) = 1$), then $ \mu$-a.e.\ $\omega_1 \in \Omega$ lies in $ M_\mu $, and the support $\Omega_\mu $ lies in the topological closure of $ M_\mu $. In addition, for many dynamical systems encountered in applications, $M_\mu$ can be a significantly ``larger'' set than $ \Omega_\mu $. In particular, for systems possessing physical measures \cite{Young02}, $M_\mu $ has positive measure with respect to an ambient measure on $ \Omega $ (e.g., Lebesgue measure), whereas $\Omega_\mu$ oftentimes has zero ambient measure (e.g., if $\Omega_\mu$ is an attractor developing under dissipative dynamics). In such cases, the methods will converge from an experimentally accessible set of initial conditions that can lie outside of $\Omega_\mu$. 
\paragraph*{Uniform convergence on $\mathcal{U}$} We assume throughout that the basic assumptions stated above hold. For simplicity, we will assume that for the given hypothesis space dimension $\ell$, all eigenvalues $\lambda_1, \ldots, \lambda_{\ell}$ are simple (if this is not the case, the argument presented below can be modified using appropriate projector operators onto eigenspaces of $G$ and $G_n$). 

Since the ideal target function $ f_{\tau,\ell} $ from~\eqref{eq:firstKPCR} and the empirical target function~\eqref{eqKAF} are linear combinations of $\ell < \infty$ continuous functions, $ \psi_i $ and $\psi_{i,n}$, respectively, convergence of $f_{\tau,\ell,n}$ to $f_{\ell}$ in $C(\mathcal{U})$ norm will follow if it can be shown that, as $n \to \infty$ and for each $i \in \{ 1, \ldots, \ell\} $, (i) the eigenvalues $\lambda_{i,n} $ converge to $\lambda_i$;  (ii) the RKHS functions $\psi_{i,n}$ converge, up to multiplication by a constant phase factor, to $\psi_i$ in $C(\mathcal{U})$ norm; and (iii) each of the expansion coefficients $\alpha_{i,n}(\tau)$ converges to $\alpha_i(\tau)$. The first two of these claims are a consequence of the following lemma, which is based on \citep[][Theorem~15]{VonLuxburgEtAl08}, \citep[][Corollary~2]{DasGiannakis19}, and \citep[][Theorem~7]{GiannakisEtAl19}. 

\begin{lem}{\label{lemSpecConv} Under the basic assumptions for convergence, the following hold:
        \begin{enumerate}[(i),wide]
            \item For each nonzero eigenvalue $\lambda_i$ of $G$, $\lambda_{i,n}$ converges to $\lambda_i$ as $n \to \infty$.
            \item For every RKHS function $\psi_i$ corresponding to $\lambda_i >0$, there exist complex numbers $c_{i,n}$ of unit modulus, such that $\lim_{n\to\infty} \lVert \psi_{i,n} - c_{i,n} \psi_i \rVert_{C(\mathcal{U})} = 0$.
\end{enumerate}}
\end{lem}

\begin{rk} \label{rkRosasco} In \cite{RosascoEtAl10}, Rosasco et al.\ approach the problem of establishing spectral convergence of the empirical integral operators $G_n = \iota_n \iota_n^*$ associated with a fixed ($n$-independent) kernel $k$ by considering the operators $\tilde G = \iota_n^* \iota_n $ acting on the corresponding RKHS, $\mathcal K$. They show that for i.i.d.\ training data sampled from $\mu_X$, as $n \to \infty$ these operators converge in Hilbert-Schmidt norm, and thus in spectrum, to the integral operator $ \tilde G = \iota^* \iota$, and provide explicit rates of convergence. Aside from loosing the precise error bounds that the i.i.d.\ assumption affords, the weaker ergodicity assumption employed in this work could be used to establish convergence in Hilbert-Schmidt norm of the data-driven operators associated with the sampling measures $\mu_{X,n}$ along orbits of the dynamics, for which the data are not independent. An advantage of this approach is that one does not need to introduce $C(\mathcal U)$ as an auxiliary comparison space for the operators $G_n$ and $G$ (which act on different spaces). Moreover, convergence of $\tilde G_n$ in Hilbert-Schmidt norm is stronger than the type of operator convergence considered in  \cite{VonLuxburgEtAl08,DasGiannakis19,GiannakisEtAl19} (namely, collective compact convergence), which leads to Lemma~\ref{lemSpecConv}. At the same time, however, a potential limitation of the Hilbert-Schmidt approach is that it requires the existence of a fixed RKHS, and as previously discussed, in many cases  the kernels $k_n$ and corresponding RKHSs  $\mathcal K_n$ depend on the training data. For this reason, we have opted to work in the more general setting of Lemma~\ref{lemSpecConv} despite a somewhat weaker convergence result, but we should point out that the results of \cite{RosascoEtAl10} are available as an option when KAF is implemented with data-independent kernels.    
\end{rk}

Next, let $\tilde T_n : K_n \to \mathcal{K}_n $,  be the empirical Nystr\"om extension operator on $ K_n := \iota_n \mathcal{K}_n$, defined analogously to $\tilde T $ from Section~\ref{secHypothesis}. Also, for any probability measure $\rho : \mathcal{F} \to [ 0, 1 ] $, let $ \rho( f ) = \int_\mathcal{X} f \, d\rho $, where $ f \in \mathbb{L}^1(\rho)$. To verify convergence of the expansion coefficients $\alpha_{i,n}(\tau)$, note that Lemma~\ref{lemSpecConv}(ii) implies that for each $ i $ such that $\lambda_i > 0 $,  the continuous representatives of $\phi_{i,n}$, given by $ \varphi_{i,n} = \tilde T_n \phi_{i,n} = \psi_{i,n} / \lambda_{i,n}^{1/2} $, converge in $C(\mathcal{U} )$ norm and up to phase to the continuous representative   $\varphi_{i}= \tilde T \phi_i = \psi_{i} / \lambda_i^{1/2}$ of $ \phi_i$. Moreover, because the products $\alpha_{i,n}(\tau) \psi_{i,n}$ are invariant under multiplication of $\phi_{i,n}$ by a constant phase factor, without loss of generality, we may assume that the $ c_{i,n}$ in Lemma~\ref{lemSpecConv} are all equal to 1. Then, for any $ \tau = q \, \Delta t$ with $ q \in \mathbb N_0 $, we have
\begin{align*}
    \alpha_{i,n}(\tau) &= \mu_n( (\phi^*_{i,n} \circ X ) ( U^{q}_n \iota_n Y) ) \\
    &= \frac{1}{n} \sum_{j=1}^{n-q} \varphi^*_{i,n}(x_j) U^\tau Y(\omega_{j}) \\  
    &= \frac{n-q}{n} \mu_{n-q}( ( \varphi^*_{i,n} \circ X ) U^\tau Y ),
\end{align*}
and defining $\tilde \mu_n = (n'/n hypothesis space dimension) \mu_n$, it follows that
\begin{multline*}
    \lvert \alpha_{i,n}(\tau) - \alpha_i(\tau) \rvert  \\
    \begin{aligned}
        & = \lvert \mu_n( ( \phi_{i,n}  \circ X ) U^{q}_n \iota_n Y ) - \mu((\phi_i \circ X) U^\tau Y ) \rvert \\
        &= \lvert  \tilde\mu_n( ( \varphi_{i,n} \circ X ) U^\tau Y ) - \mu( (\varphi_i \circ X) U^\tau Y ) \rvert \\
        & \leq \lvert \tilde \mu_n( [ ( \varphi_{i,n} - \varphi_i ) \circ X]  U^\tau Y )  \\
        & \quad + \lvert ( \tilde \mu_n - \mu )((\varphi_i \circ X ) U^\tau Y) \rvert, \\
        & \leq  C_Y \lVert \varphi_{i,n} - \varphi_i \rVert_{C(\mathcal{U})} \\ 
         & \quad + \lvert ( \tilde \mu_n - \mu )( ( \varphi_i \circ X ) U^\tau Y) \rvert.
     \end{aligned}
\end{multline*}
In the last line above, the first term converges to 0 by uniform convergence  of $ \varphi_{i,n} $ to $ \varphi_{i}$ on $\mathcal{U} $, and the second term by weak convergence of $\mu_{n}$ to $ \mu$, so we conclude that $\alpha_{i,n}(\tau)$ converges to $\alpha_i(\tau)$. Moreover, by continuity of the dynamics and covariate and response variables, the convergence is uniform with respect to $\tau $ lying in compact sets.

We summarize the main results of Sections~\ref{secGenError} and~\ref{secSampleError} in the following theorem:
\begin{thm}
    \label{thmEmpiricalConv}
    Under the basic assumptions for convergence, for every lead time $ \tau = q \, \Delta t $, $ q \in \mathbb N_0$, and hypothesis space dimension $\ell$ such that $\lambda_{\ell} >0$, the KAF target function $f_{\tau,\ell,n} \in \mathcal{K}_n$ converges as $n \to \infty$ to the ideal target function $f_{\tau,\ell} \in \mathcal{K}$, uniformly on $\mathcal{U}$. Moreover, if the reproducing kernel $k$ of $ \mathcal{K} $ is $L^2(\mu_X)$-strictly-positive-definite, then by Theorem~\ref{thmConv}, $f_{\tau,\ell,n}$ converges to the regression function $Z_\tau $ associated with the conditional expectation, $ \mathbb E[ U^\tau Y \mid X] = Z_\tau \circ X $, in the sense of the iterated limit
    \begin{displaymath}
        \lim_{\ell\to\infty}\lim_{n\to\infty} f_{\tau,\ell,n}  = \lim_{\ell\to\infty} f_{\tau,\ell} = Z_\tau.
    \end{displaymath}
    Here, the $n\to\infty$ and $\ell \to \infty$ limits are taken in $C(\mathcal{U})$ and $L^2(\mu_X)$ norm, respectively. Moreover, the convergence is uniform with respect to $ \tau $ lying in compact sets.  
\end{thm}

\subsection{\label{secMixing}Mixing and loss of predictability}

Before closing Section~\ref{secError}, we discuss some aspects of the long-time behavior of the conditional expectation and the KAF target functions in the presence of mixing dynamics, which will be useful in our interpretation of the L63 experiments in Section~\ref{secL63}. First, we recall that the measure-theoretic definition of mixing \citep[e.g.,][]{Walters81} can be equivalently stated as the condition that  for any $g,h \in L^2(\mu)$, 
\begin{displaymath}
    \lim_{\tau \to \infty} \langle U^{\tau*} g, h \rangle_{L^2(\mu)} = \left( \int_{\Omega} g^* \, d\mu \right) \left( \int_\Omega h \, d\mu \right).
\end{displaymath}
Thus, under mixing dynamics, inner products of the form $ \langle U^{\tau*} g, h \rangle_{L^2(\mu)} $, which can be thought of as temporal cross-correlation functions, converge to constants equal to products of the expectation values $ \mathbb{E}[g^*]= \int_\Omega g^* \, d\mu $ and $ \mathbb{E}[h] = \int_\Omega h \, d\mu $. Using the projection representation of the conditional expectation in~\eqref{eqEProj}, it then follows that
\begin{align*}
    \lim_{\tau\to\infty} \langle g, \mathbb{E}[U^\tau Y \mid X] \rangle_{L^2(\mu)} &= \lim_{\tau\to\infty}\langle g, \Pi_X U^\tau Y \rangle_{L^2(\mu)}\\
    &= \lim_{\tau\to\infty} \langle U^{\tau*} \Pi_X g, Y \rangle_{L^2(\mu)} \\
    &  = \mathbb{E}[ ( \Pi_X g)^* ] \mathbb{E}[  Y ].
\end{align*}
Therefore, because $g $  in the above is arbitrary, and $\Pi_X$ leaves constant functions invariant, we conclude that $\mathbb{E}[U^\tau Y \mid X]$ converges weakly to a constant function equal to $ \mathbb{E}[Y]$, i.e., 
\begin{multline*}
    \lim_{\tau\to\infty} \langle g, \mathbb{E}[U^\tau Y \mid X] - \mathbb{E}[Y] 1_\Omega \rangle_{L^2(\mu)} \\= \lim_{\tau\to\infty} \langle g, \Pi_X ( U^\tau Y - \mathbb{E}[Y] 1_\Omega) \rangle_{L^2(\mu)} = 0,
\end{multline*}
where $ 1_\Omega $ is the function on $ \Omega$ equal everywhere to 1. We interpret this behavior as a loss of predictability due to mixing dynamics. 

Observe now that the $L^2(\mu)$ element $ f_{\tau,\ell} \circ X $, where $ f_{\tau,\ell} $ is the ideal target function, can be expressed as $ \Pi_{\tilde H_l} U^{\tau} Y $, where $ \Pi_{\tilde H_l } $ is the orthogonal projection on $L^2(\mu) $ mapping into the pullback $ \tilde H_\ell = H_\ell \circ X $ of the $\ell$-dimensional hypothesis space $H_\ell$ into $L^2_X(\mu)$. If $ \tilde H_\ell $ contains constant functions, then it follows from similar arguments as above, in conjunction with the fact that $ \tilde H_\ell $ is finite-dimensional, that as $\tau \to \infty $, $f_{\tau,\ell} $ converges in $L^2(\mu)$ norm (and not merely weakly) to $\mathbb{E}[Y]$. We will discuss practical ways of ensuring that $ \tilde H_\ell $ always contains constant functions, ensuring in turn this type of long-term statistical consistency with the infinite-dimensional case, using Markov-normalized reproducing kernels in Section~\ref{secKernel}.    

With regards now to the empirical target function, since the $n \to \infty$ convergence of $f_{\tau,\ell,n}$ to $ f_{\tau,\ell} $ may not be uniform with respect to $\tau \in \mathbb{R} $, we cannot use this result to make a statement about the relation between $f_{\tau,\ell,n} $ and  $\mathbb{E}[Y]$ as $ \tau \to \infty $. Nevertheless, it is still possible to ensure (through Markov normalization of the kernel) that, at fixed $n $,  $ f_{\tau,\ell,n} \circ X $ lies in a finite-dimensional subspace of $L^2(\mu) $ containing constant functions. In that case, for large-enough $ n $, and long-enough, but bounded, $ \tau $, we can expect $ f_{\tau,\ell,n} $ to be an approximately constant function equal to $ \mathbb{E}[Y] $.   

\begin{rk}
    Time series prediction techniques can generally be categorized as being \emph{direct} or \emph{iterated} methods \cite{Sauer92}. The KAF target function $f_{\tau,\ell,n}$ from~\eqref{eqKAF} provides direct prediction, in the sense that every lead time $\tau $ has a distinct associated forecast function, which is evaluated once at the given initial data to yield a prediction. In contrast, in iterated prediction, one sets a timestep $ \Delta t > 0 $, and constructs a function $ g : \mathcal X \to \mathcal X $ that propagates the covariate over that timestep. Then, to obtain prediction with lead time $ \tau = q \, \Delta t$, $q \in \mathbb N$, the function $ g $ is iteratively evaluated $q$ times, $ x_{j+1} = g( x_j )$, using the covariate $x \in \mathcal X $ observed at forecast initialization as the initial condition $ x_0 = x $, and the result is fed into a function $ f : \mathcal X \to \mathcal Y$, yielding a forecast $ y = f( x_q ) $ of the response. Direct methods are more general than iterated methods, as the latter are generally based on some type of Markovianity assumption for the training data. Indeed, in Theorem~\ref{thmEmpiricalConv}, the convergence of $f_{\tau,\ell,n}$ to the $L^2$-optimal conditional expectation was established without invoking any assumption about the dynamics in covariate space. In contrast, iterated methods effectively construct a surrogate dynamical model on $\mathcal X$, which if successful, may provide access to long-term statistics in ways that are not possible by direct methods. For example, the time series $ (q \, \Delta t) \mapsto ( f \circ g^q)(x)$ produced by an iterated method need not converge to a constant as $q\to\infty$, but may exhibit variance and higher-order statistical moments resembling those induced by the invariant measure. See \cite{SmallJudd98} for examples of iterated models using radial basis functions.
\end{rk}

\section{\label{secExt}Extensions}

This section shows how the KAF/KPCR learning framework presented thus far can shed light on other aspects of the kernel approach other than using leading principal components (eigenfunctions) to approximate the conditional expectation of observables. The first extension (Section~\ref{secKRR}) describes how KRR may be characterized as resulting from the same variational problem as that of KPCR, albeit with a nonlinear, rather than linear, hypothesis space. The second extension (Section~\ref{secAsym}) shows how KAF can be implemented using a class of non-symmetric kernels. The third extension (Section~\ref{secProb}) shows how quantities other than the conditional expectation, such  as the conditional probability and estimates of the forecast error, may also be approximated, and what their utility may be in practical problems.

\subsection{\label{secKRR}Kernel ridge regression (KRR)}
In KRR \cite{saunders1998ridge}, the hypothesis space is a closed ball of radius $R$ in the RKHS $\mathcal{K}$;
\begin{displaymath}
\mH_R = \left\{ f \in \mathcal{K} : \lVert f\rVert_\mathcal{K} \le R\right\}.
\end{displaymath}
Note that $\mH_R$ is not a linear subspace of $\mathcal{K}$, and thus the projection $\Pi_{H_R} : L^2(\mu_X) \to L^2(\mu_X)$ mapping into $ H_R = \iota \mathcal{H}_R$ is a nonlinear operator. Although representations for this particular nonlinear operator are known \cite{cucker2007learning}, those for the inverse $\iota_{\mH_R}^{-1}$ are generally intractable. However, by using Lagrange multipliers, optimization over $\mH_R$ may be transformed into a linear problem. In particular, the problem of minimizing $\lVert \iota f - Z_\tau\rVert_{L^2(\mu_X)}$ such that $\lVert f \rVert_\mathcal{K} \le R$ is a constrained optimization problem for which there exists a parameter $\eta_R > 0$, dependent on $R$, such that the penalized optimization problem
\begin{displaymath}
\min_{f\in\mathcal{K}} \rVert^2 \iota f - Z_\tau \rVert^2_{L^2(\mu_X)} + \eta_R  \lVert f \rVert_\mathcal{K}^2
\end{displaymath}
is an equivalent formulation. The solution to this problem is known to be \cite{cucker2007learning}
\begin{displaymath}
f_{\tau,R} = \iota^* (G + \eta_R \Id)^{-1} Z_\tau.
\end{displaymath}
The empirical solution $f_{\tau,R,n} \in \mathcal{K}_n$, meanwhile, is given by
\begin{equation}
    \label{eqKRR}
f_{\tau,R, n}(x) = \bk(x)^* (\bm G + \eta_R {\bf I})^{-1} \bm y_\tau,
\end{equation}
where $ \bm G $, $ \bm k $, and $ \bm y_\tau $ are the kernel matrix, kernel vector, and analog vector from Section~\ref{secTarget}. As with KPCR, the KRR target function $f_{\tau,R}$ also converges in mean square to the conditional expectation, in the sense that, as the regularization parameter $\eta_R$ is decreased to zero, $\lVert \iota f_{\tau,R} - Z_\tau\rVert_{L^2(\mu_X)} $ converges to 0 if the kernel $k$ is $L^2(\mu_X)$-strictly-positive (cf.\ Theorem~\ref{thmConv}). Moreover, under the assumptions stated in Section~\ref{secSampleError}, $f_{\tau,R,n}$ converges to $f_{\tau,R}$ in $C(\mathcal{U}) $ norm, as $ n \to \infty $, so that an analog of Theorem~\ref{thmEmpiricalConv} holds for $f_{\tau,R,n}$.

Note that, unlike KPCR, the KRR estimator in~\eqref{eqKRR} does not require eigendecomposition of $\bm G $, and only depends on kernel values (thus, it can be thought of as making direct use of the ``kernel trick''; see Section~\ref{secRKHS}.) Still, the standard implementation of KRR relies on a computationally expensive full inversion of the kernel matrix, whose eigenvalues are perturbed away from zero by some regularizing parameter $\eta_R$. A hybrid approach is to employ a low-rank approximation as in KPCR, while perturbing the eigenvalues away from zero as in KRR. With the notation of~\eqref{eqAMat}, this leads to the target function
    \begin{displaymath}
        f_{\tau,\ell,R,n}(x) = \bk(x)^T \bm \Phi_{\ell} ( \bm \Lambda_{\ell} + \eta_R \bm I )^{-1} \bm \Phi_{\ell}^* / n^2.
    \end{displaymath}
All of the KPCR, KRR, and hybrid estimators approximate the conditional expectation when the parameters are sufficiently relaxed, but the rates of convergence may differ. In general, KRR is useful when insensitivity to noise is desired, but it can be computationally expensive as it involves full matrix inversion. KPCR, on the other hand, can converge very rapidly when it turns out that the regression function lies in the leading eigenspaces of $G$.

\subsection{\label{secAsym}Non-symmetric kernels}

The KAF formulation presented in Sections~\ref{secKAF} and~\ref{secError} makes heavy use of RKHSs, and is therefore restricted to Hermitian, positive-definite kernels. Yet, many popular kernel-based algorithms utilize non-symmetric kernels, typically constructed by normalization of symmetric kernels. Examples include normalized graph Laplacians \cite{VonLuxburgEtAl08} and Markov kernels approximating heat kernels on manifolds \cite{HeinEtAl05,CoifmanLafon06,BerryHarlim16,BerrySauer16b}. We now describe an extension of KAF to a class of non-symmetric kernels, whose corresponding integral operators on $L^2(\mu)$ are related to integral operators associated with symmetric positive-definite kernels by similarity transformations. We will see that target functions can still be constructed using Nystr\"om-type extensions into RKHSs, with extension operators derived from integral operators induced by non-symmetric kernels.
    
Specifically, we let $ w : \mathcal X \times \mathcal X \to \mathbb R $ be a continuous, positive-definite bivariate function on covariate space, not necessarily symmetric, such that 
\begin{equation}
    \label{eqBalance}
    d(x) w(x,x') = d(x')w(x'x), \quad \forall x,x'\in \mathcal X,
\end{equation}
for a continuous, strictly positive function $ d : \mathcal X \to \mathbb R $. Notice the similarity between~\eqref{eqBalance} and the detailed balance relation in reversible Markov chains. We have: 

\begin{prop}\label{propAsym} With notation as above, and if~\eqref{eqBalance} holds, the function $k : \mathcal X \times \mathcal X \to \mathbb R $ with
    \begin{displaymath}
        k(x,x') = w(x,x') / d(x')
    \end{displaymath}
    is a continuous, symmetric, positive-definite kernel. Moreover, denoting the corresponding RKHS, inclusion map, and Nystr\"om operator by $\mathcal K$, $\iota : \mathcal K \to L^2(\rho)$, $ \tilde T = \iota^+$, respectively, where $\rho$ is any compactly supported finite Borel measure on $\mathcal X$, the following hold.
    \begin{enumerate}[(i),wide]
    \item The integral operator 
        \begin{displaymath}
             W : f \mapsto \int_{\mathcal X} w( \cdot, x ) f(x) d\rho(x )   
        \end{displaymath}
        is well-defined as a bounded operator $W : L^2(\rho) \to \mathcal K$.
    \item The integral operator $ J = \iota W $ is a trace-class operator on $L^2(\mu_X)$, with real eigenvalues $ \eta_1 \geq \eta_2 \geq \cdots \searrow 0^+$. Moreover, there exists a Riesz basis $ \{ \xi_1, \xi_2, \ldots, \} $ of $L^2(\mu_X) $ and a dual basis $ \{ \xi'_1, \xi'_2,\ldots \}$ with $ \langle \xi'_i, \xi_j \rangle_{L^2(\mu_X)} = \delta_{ij}$ consisting of eigenfunctions of $ J $ and $ J^*$, respectively, i.e., 
        \begin{displaymath}
            J \xi_j = \eta_j \xi_j, \quad J^* \xi'_j = \eta_j \xi'_j.
        \end{displaymath}
    \item The domain $D(J^+)$ of the Moore-Penrose pseudoinverse of $J$ is a dense subspace of $D(\tilde T)$, and on this subspace the Nystr\"om operator $\tilde T$ takes the form $ \tilde T\rvert_{D(J^+)} = W J^+ $. As a result, for every $f \in D(\tilde T)$ we have
        \begin{equation}
            \label{eqNystAlt}
            \tilde T f = \sum_{i:\eta_i > 0} \frac{\langle \xi'_i, f \rangle_{L^2(\mu_X)}}{\eta_i^{1/2} }\vartheta_i, 
        \end{equation}
        where $ \vartheta_i = W \xi_i / \eta_i^{1/2}$ are orthogonal functions in $\mathcal K$. 
\end{enumerate}
\end{prop}

Proposition~\ref{propAsym} was inspired by \cite[][Section~4]{RosascoEtAl10}, where an auxiliary RKHS analogous to $\mathcal K$ was used to establish spectral convergence for a class of non-symmetric graph Laplacian operators. Here, our perspective is somewhat different as we introduce $\mathcal K $  as a consequence of the detailed balance condition in~\eqref{eqBalance}, rather than assuming its existence a priori and deducing from it a non-symmetric kernel such as $w$ (as done in~\cite{RosascoEtAl10}). Moreover, our objective here is not to establish spectral convergence of integral operators acting on $\mathcal K$ (since our kernels are typically data-dependent; see Remark~\ref{rkRosasco}), but rather to identify an RKHS whose corresponding Nystr\"om operator is computable using integral operators associated with non-symmetric kernels.  Indeed, similarly to \cite{RosascoEtAl10}, a key aspect of Proposition~\ref{propAsym} is that the Nystr\"om operator for $\mathcal K$ is constructed using integral operators associated with $w$ and their pseudoinverses, without invoking the operator $\iota^*$ associated with the reproducing kernel $k $ of $ \mathcal K$. 

\begin{proof} 
    See \ref{appAsym}.
\end{proof}

The eigenvalues $\eta_i$ appearing in Proposition~\ref{propAsym} are different from the eigenvalues $\lambda_i $ of the integral operator $ G = \iota \iota^*$ associated with $k$. Instead, as can be directly verified, the $\eta_i$ coincide with the eigenvalues of a self-adjoint integral operator $\hat G$ on $L^2(\rho)$ associated with a different reproducing kernel than $k$, to which $J $ is related by a similarity transformation. Specifically, observe that under detailed balance, the kernel $\hat k : \mathcal X \times \mathcal X \to \mathbb R$ given by 
\begin{equation}
    \label{eqKHat}
    \begin{aligned}
        \hat k(x,x') &= d^{1/2}(x) w(x,x') d^{-1/2}(x') \\
        &= d^{1/2}(x)k(x,x')d^{1/2}(x')
    \end{aligned}
\end{equation}
is also positive-definite, and has an associated RKHS $\hat{\mathcal K}$ with restriction operator $ \hat \iota : \hat{\mathcal K} \to L^2(\rho)$ and self-adjoint integral operator $ \hat G = \hat \iota \hat \iota^* : L^2(\rho) \to L^2(\rho)$. Letting $M_d : L^2(\rho) \to L^2(\rho) $ be the bounded multiplication operator by the continuous function $d$, i.e., $ M_d f = d f $, we can express $J$ as a similarity transformation of $\hat G $, namely 
\begin{equation}
    \label{eqJGHat}
    J = M_d^{-1/2} \hat G M_d^{1/2}, 
\end{equation}
where $M_d^{1/2}$ and $M_d^{-1/2}$ are both bounded multiplication operators by the strict positivity and continuity of $ d $ and compactness of the support of $\mu_X$. It then follows that $J$ and $G$ have the same eigenvalues $\eta_i$, and we can construct the biorthonormal bases $ \xi_i $ and $ \xi_i' $ by first computing orthonormal eigenvectors $ \hat G_i \hat \phi_i = \eta_i \hat \phi_i$ (taking advantage of specialized solvers for self-adjoint operators), and setting
\begin{equation}
    \label{eqXi}
    \xi_i = d^{-1/2} \hat \phi_i, \quad \xi_i' = d^{1/2} \hat \phi_i. 
\end{equation}
Similarly, we have 
\begin{equation}
    \label{eqTheta}
    \vartheta_i = d^{-1/2} \hat \psi_i, \quad \hat \psi_i = \hat \iota^* \hat \phi_i / \eta_i^{1/2}.
\end{equation}
Note that $G = M_d^{-1/2} \hat G M_d^{-1/2}$, so apart from special cases, $G$ and, $\hat G$ are not similar operators.

As in the case of symmetric kernels, when forecasting with the class of non-symmetric kernels considered in this section we construct target functions analogously to~\eqref{eq:theFormula} and~\eqref{eqKAF} by truncating the expansion in~\eqref{eqNystAlt} to $\ell$ eigenfunctions. 

\subsection{\label{secProb}Conditional variance and conditional probability}

In forecasting applications, it is important to be able to perform uncertainty quantification; that is, estimate the error of the target function. Moreover, besides point forecasts of a given response variable, it is oftentimes of interest to predict the probability of occurrence of events defined in terms of the response meeting certain criteria (e.g., exceeding a specified threshold). We now discuss how KAF techniques can be employed to carry out these tasks.

First, regarding error estimation, consider the $L^2(\mu)$ observable
\begin{displaymath}
    \beta_\tau = \lvert U^\tau Y - \mathbb{E}[U^\tau Y \mid X] \rvert^2,
\end{displaymath}
which measures the square error of the conditional expectation (and thus, the optimal target function $ Z_\tau$).  The conditional expectation of $ \beta_\tau$ with respect to $X$ is known as the conditional variance,
\begin{displaymath}
    \var[ U^\tau Y \mid X] = \mathbb{E}[\beta_\tau \mid X], 
\end{displaymath}
and satisfies
\begin{displaymath}
    \int_\Omega \var[ U^\tau Y \mid X ] \, d\mu = \sigma_\tau
\end{displaymath}
by construction. Thus, the conditional variance is equal in expectation to the intrinsic error from~\eqref{eqSigmaRho}, providing an unbiased estimator of the square forecast error. Moreover, being a conditional expectation, $ \var[ U^\tau \mid X ]$ is expressible as the pullback of a unique element $ W_\tau \in L^2(\mu_X)$, such that $ \var[ U^\tau \mid X ] = W_\tau \circ X$, which can be empirically approximated using KAF. Specifically, applying KAF to the function $ \beta_{\tau,\ell,n} = \lvert U^\tau Y - f_{\tau,\ell,n} \circ X \rvert^2 $  leads to the estimator $ s_{\tau,\ell,n} \in \mathcal{H}_{\ell,n} $ of $W_\tau$ given by 
\begin{displaymath}
    s_{\tau,\ell,n} = \bm k(x)^* \bm A_\ell \bm \beta_{\tau,\ell,n}, 
\end{displaymath}
where $ \bm \beta_{\tau,\ell,n } = [ \beta_{\tau,\ell,n}(\omega_1), \ldots, \beta_{\tau,\ell,n}(\omega_n) ]^T $ is a column vector in $\mathbb{C}^n $ containing the values of $ \beta_{\tau,\ell,n} $ on the training states $ \omega_i $ (cf.\ the analog vector $ \bm y_\tau  $ in~\eqref{eqAnalogVec}). Because, as with many projection methods, $ s_{\tau,\ell,n}(x) $ is not guaranteed to be non-negative, in practice we perform error estimation using
\begin{equation}
    \label{eqErrEst}
    \varepsilon_{\tau,\ell,n}(x) = \lvert s_{\tau,\ell,n}(x) \rvert^{1/2}.
\end{equation}
The function $ \varepsilon_{\tau,\ell,n} \circ X $ then converges in the limit of large data to $ \sqrt{ \var[ U^\tau Y \mid X ] } $, analogously to the convergence of $ f_{\tau,n,\ell} \circ X $ to $ \mathbb{E}[U^\tau Y \mid X] $ in Theorem~\ref{thmEmpiricalConv}. More generally, note that for any function $\Gamma : \mathbb{C} \to \mathbb{C} $, such that $ \Gamma \circ Y $ lies in $\mathbb{L}^2(\mu)$, the  conditional expectation $\BE[ U^\tau( \Gamma \circ Y ) \mid \bm X]$ is approximated by
\begin{equation}
    \label{eqKAF_G}
    g_{\tau,\ell,n} = \bm k(x)^* \bm A_\ell \Gamma(\bm y_\tau),
\end{equation}
where $\Gamma(\bm y_\tau ) $ is the column vector in $ \mathbb{C}^n$ obtained by element-wise application of $ \Gamma $ to the analog vector $ \bm y_\tau $.  

Next, turning to approximations for conditional probability, let $ \Theta \in \mathcal{F} $ be an event (i.e., a measurable subset of $\Omega$), defined through certain conditions on $Y(\omega) $ being met. For instance, in the forecasting of rare or extreme events, one might employ a formulation such as  
\begin{equation}
    \label{eqSThresh}
    \Theta = \{ \omega \in \Omega: Y(\omega) > \theta \},
\end{equation}
where $\theta$ is a large threshold parameter. 

Every event $\Theta$ has an associated indicator function  $ \chi_\Theta \in \mathbb{L}^2(\mu)$, evolving under the action of the Koopman operator as 
\begin{displaymath}
    U^\tau \chi_\Theta = \chi_{\Theta_\tau}, \quad \Theta_\tau = A^{-\tau}(\Theta).
\end{displaymath}
Note, in particular, that every point lying initially in $ \Theta_\tau $ will be mapped into $\Theta$ after dynamical evolution over time $ \tau $.  The conditional expectation  
\begin{displaymath}
\mathbb{P}[ \Theta_\tau \mid X ] := \mathbb{E}[ U^\tau \chi_\Theta \mid X ]
\end{displaymath}
then gives the conditional probability for $ \Theta $ to occur at lead time $ \tau $ given $X$. In the context of KAF, approximations for conditional probability are obtained by setting $\Gamma$ in~\eqref{eqKAF_G} to be the indicator function of the set $ Y(\Theta) \subset \mathbb{C} $,  leading to the target function 
\begin{displaymath}
    \tilde g_{\tau,\ell,n} = \bm k (x)^T \bm A_\ell \chi_{Y(\Theta)}(\bm y_\tau). 
\end{displaymath}
Because $\tilde g_{\tau,\ell,n}$ is not guaranteed to take values in the interval $ [ 0, 1] $, in order to obtain meaningful forecasts of conditional probability we threshold it, leading to the estimator
\begin{equation}
    \label{eqKAFProb}
    g_{\tau,\ell,n}(x) = 
    \begin{cases}
        1, & g_{\tau,\ell,n}(x) > 1, \\
        g_{\tau,\ell,n}(x), & 0 < g_{\tau,\ell,n}(x) \leq 1, \\
        0, & g_{\tau,\ell,n}(x) \leq 0,
    \end{cases}
\end{equation}
where $ g_{\tau,\ell,n} \circ X $ approximates $ \mathbb{P}[\Theta_\tau \mid X]$ analogously to Theorem~\ref{thmEmpiricalConv}. For example, for the event in~\eqref{eqSThresh}, $ g_{\tau,\ell,n}(x) $ estimates the likelihood that $ Y $ will exceed $ \theta $ at lead time $\tau $, given the covariate value $x \in \mathcal{X}$.

\section{\label{secKernel}Choice of kernel}

In this section, we discuss practical guidelines for choosing the kernel $k_n $ on covariate space $ \mathcal{X}$ employed in KAF. 

\subsection{\label{secPD}Strictly positive-definite kernels} 

As a general guideline, in order to ensure that the empirical target function $f_{\tau,\ell,n}$ from Definition~\ref{defnEmpTarget} converges to the regression function $Z_\tau $ from Definition~\ref{defnRegression}  for an arbitrary response variable $Y$ and lead time $ \tau $ (i.e., Theorem~\ref{thmEmpiricalConv} holds), the empirical kernels $k_n$ should converge, as $n \to \infty$, to an $L^2(\mu_X)$-strictly-positive kernel $k$, uniformly on the compact set $\mathcal{U} \subseteq \mathcal{X}$. Because every Mercer kernel which is strictly positive-definite on the support of a compactly supported Borel probability measure $ \rho $  is $L^2(\rho)$-strictly-positive (see Section~\ref{secRKHS}), a convenient way of ensuring $L^2(\mu_X)$-strict-positivity of $k $ is to work with empirical kernels $ k_n $ whose restrictions on  $ \supp \mu_{X,n} =  \{ x_1, \ldots, x_n \} $ are strictly positive-definite for every $ n \in \mathbb{N} $.  For example, in the case $\mathcal{X} = \mathbb{R}^m$, it is known that radial Gaussian kernels are strictly positive-definite on the whole of $\mathcal{X}$ \cite{Micchelli86}. Therefore, one can work with
\begin{equation}
    k_n(x,x') = k(x,x' ) = e^{- \lVert x - x' \rVert^2 / \epsilon}
    \label{eqKGauss}
\end{equation}
for some positive bandwidth parameter $ \epsilon $, and the conditions of Theorem~\ref{thmEmpiricalConv} will be satisfied. 

The radial Gaussian kernel in~\eqref{eqKGauss} will be employed in the circle example in Section~\ref{secCircle}. It is an instance of a local kernel \cite{BerrySauer16b} of the form
\begin{equation}
    \label{eqKLocal}
    k(x,x') = h(L(x,x') /\epsilon),
\end{equation}
where $L $ is a continuous, positive, symmetric function on $\mathcal X \times \mathcal X \to \mathbb R $, and $h : \mathbb R_+ \to \mathbb R_+ $ a strictly positive, bounded, continuous, shape function with rapid (faster-than-polynomial) decay at infinity. In the case of~\eqref{eqKGauss}, we have $L(x,x') = \lVert x - x' \rVert^2$ and $h(u) = e^{-u} $. See Ref.~\cite{Genton01} for additional examples of kernels commonly used in machine learning applications.          

\begin{rk}On the finite-dimensional linear covariate space $\mathcal X = \mathbb{R}^m$, the covariance kernel, $k(x,x') = x \cdot x'$, which is employed in the proper orthogonal decomposition \cite{Kosambi43,HolmesEtAl96} and linear inverse modeling techniques \cite{Penland89}, is not $L^2(\mu_X)$-strictly-positive. Indeed, one can verify that for this choice of kernel, the corresponding integral operators $G_n$ and $G $ are of at most rank $m$ \cite[][Section~9]{DasGiannakis18}, thereby bounding the dimension of the hypothesis spaces $\mathcal{H}_{\ell,n}$ and $ \mathcal{H}_{\ell} $ by $ m $. Thus, unless $Z_\tau$ happens to lie in the span of the leading $d $ eigenfunctions of $G $ (which are, in this case, linear functions on $\mathcal{X}$), the empirical target function $ f_{\tau,\ell,n} $ will fail to converge to $Z_\tau $. 
\end{rk}

\subsection{\label{secMarkov}Variable-bandwidth, Markov kernels}

Next, we discuss two modifications of the radial Gaussian kernel on $\mathbb{R}^d$, which can play a fairly substantial role in improving the robustness of the hypothesis space, particularly for data with strong contrasts in sampling density (e.g., the L63 example in Section~\ref{secL63}). 

\paragraph*{Variable bandwidth} Our first modification is to introduce a strictly-positive, continuous bandwidth function $ r_n : \mathcal{X} \to \mathbb{R}$, turning~\eqref{eqKGauss} into a variable-bandwidth Gaussian kernel, viz. 
\begin{equation}
    k_n( x, x' ) = \exp\left( - \frac{ \lVert x - x' \rVert^2}{ \epsilon r_n(x) r_n(x') } \right). 
    \label{eqKVB}
\end{equation}
 Intuitively, the role of  $ r_n $ is to correct for variations in the sampling density of the data in covariate space. In particular, for a well conditioned kernel integral operator $G_n$, the number of datapoints lying within radius $O(\epsilon^{1/2}) $ balls centered at each datapoint should not exhibit significant variations across the dataset, yet, the standard radial Gaussian kernel from~\eqref{eqKGauss} has no mechanism for preventing this from happening. For appropriately chosen $ r_n $, the variable-bandwidth kernel in~\eqref{eqKVB} can, in effect, vary the radii of these balls to help improve conditioning. 
 
The different bandwidth functions proposed in the literature include near-neighbor distances \cite{ZelnikManorPerona04} and kernel density estimates \cite{BerryHarlim16}. In the numerical experiments of Section~\ref{secL63}, we will employ the latter approach, defining 
\begin{equation}
    \label{eqBandwidth}
    \begin{aligned}
        r_n( x ) &= q_n^{-1/\tilde m}( x ), \\
        q_n(x) &= \frac{ 1 }{ (\pi \tilde \epsilon )^{\tilde m / 2} } \int_\mathcal{X} e^{-\lVert x - x' \rVert^2 / \tilde \epsilon } \, d\mu_{X,n}(x').
    \end{aligned}
\end{equation} 
Here,  $ \tilde \epsilon $ a positive bandwidth parameter (different from $\epsilon $ in~\eqref{eqKVB}), and $ \tilde m $ a positive parameter approximating the dimension of the support $ \mathcal{X}_\mu $. The parameters $\epsilon $, $ \tilde \epsilon $, and $\tilde m $ are all determined from the data automatically. See \cite{BerryEtAl15,Giannakis19} for descriptions of this procedure, including pseudocode \citep[][Algorithm~1]{Giannakis19}.

It can be shown \cite{Giannakis19} that if $\mathcal{X}_\mu$ has the structure of a smooth manifold, with a Riemannian metric inherited from its embedding in $\mathcal{X} = \mathbb{R}^d$,  the bandwidth functions $r_n$ in~\eqref{eqBandwidth} induce a conformal change of metric, such that, in the new geometry, the measure $ \mu_X $ has uniform density relative to the Riemannian measure. That is, the conformal change of metric can be thought of as ``balancing out'' variations of the sampling density relative to the ambient-space metric, thus improving robustness to sampling errors. It should be noted that while here we do not assume that $\mathcal{X}_\mu$ has manifold structure, the balancing effect of the bandwidth functions is still expected to take place. 

\paragraph*{Symmetric Markov normalization} Our second modification of the radial Gaussian kernel is to normalize it to a $L^2(\mu_{X,n})$-strictly-positive Markov-ergodic kernel using the normalization procedure introduced in \cite{CoifmanHirn13}. This involves first computing the strictly positive, continuous functions
\begin{equation}
    \label{eqUV}
    \begin{aligned}
        u_n(x) &= \int_\mathcal{X} k_n(x,x') \, d\mu_{X,n}(x'), \\
        v_n(x) &= \int_{\mathcal{X}} \frac{k_n(x,x')}{u_n(x')} \, d\mu_{X,n}(x'),
    \end{aligned}
\end{equation} 
and then defining the Markov kernel $p_n : \mathcal{X} \times \mathcal{X} \to \mathbb{R}$, with 
\begin{equation}
    p_n(x,x') = \int_\mathcal{X} \frac{k_n(x,x'')k_n(x'',x')}{u_n(x)v_n(x'')u_n(x')} \, d\mu_{X,n}(x'').
    \label{eqBistoch}
\end{equation}
It can be readily verified that with this definition $ p_n $ acquires the Markov property, $ \int_\mathcal{X} p_n(x,x') \, d\mu_{X,n}(x') = 1$, for all $ x \in \mathcal{X} $. Moreover, it can be shown that if $ k_n $ is strictly positive-definite on $ \supp( \mu_{X,n}) $ then so is $ p_n $ \cite[][Lemma~12]{DasEtAl19}.  It can further be shown \cite{GiannakisEtAl19} that as $n \to \infty$, $ p_n $ converges in $C(\mathcal{U})$ norm to an $L^2(\mu_{X})$-strictly-positive Markov kernel $ p : \mathcal{X} \times \mathcal{X} \to \mathbb{R}$  (given by an analogous formula to~\eqref{eqBistoch}), so that the spectral convergence results in Lemma~\ref{lemSpecConv} hold with $ k_n $ and $ k$  replaced by $ p_n $ and $p$, respectively. 

In the context of KAF, a useful property of Markov kernels is that the associated integral operators $G_n$ and $ G $ have the top eigenvalue $\lambda_{1,n} = \lambda_1 = 1$ with a constant corresponding eigenfunction. This implies, in particular, that the corresponding RKHSs, $\mathcal{K}_n$ and $ \mathcal{K} $, respectively, always contain constant functions, and thus can naturally capture the mean of the response variable $U^\tau Y $. The eigenfunctions corresponding to $<1$ eigenvalues can then be thought of as capturing progressively finer-scale features of $U^\tau Y $, which are orthogonal to the mean. An illustration of this behavior is provided in Figures~\ref{figL63PhiFull} and~\ref{figL63PhiPartial}. 

In fact, in many ways, an RKHS $ \mathcal{K} $ with a Markov-ergodic reproducing kernel resembles a Sobolev space associated with a heat kernel on a manifold. Specifically, using the Nystr\"om extension operator, one can define a Dirichlet energy functional on the dense subspace $ K = \iota \mathcal{K} $ of $ L^2(\mu_X) $ that assigns a measure of roughness of functions analogous to the Dirichlet energy in Sobolev spaces. See \cite{DasEtAl19} for additional discussion on this topic. Appendix~B of that reference also contains pseudocode for computing the eigenfunctions $ \phi_{i,n} $ and associated RKHS functions $\psi_{i,n}$ for the kernel in~\eqref{eqBistoch}, which complements the pseudocode in Table~\ref{tablePseudocode} of this paper.

\paragraph*{Non-symmetric normalizations} While the class of symmetric Markov kernels in~\eqref{eqBistoch} is attractive due to its direct correspondence with RKHSs, in a variety of learning applications, including spectral clustering \cite{VonLuxburgEtAl08} and approximation of heat operators on manifolds \cite{CoifmanLafon06,Singer06,BerrySauer16b}, it is significantly more common to employ normalizations leading to non-symmetric kernels. As a popular example, we mention here the diffusion maps algorithm \cite{CoifmanLafon06}, which is based on the class of Markov kernels $p_n : \mathcal X \times \mathcal X \to \mathbb R$ with 
\begin{align*}
    p_n(x,x') &= \frac{\kappa (x,x')}{ v_n(x) u^\alpha_n(x')}, \\
    u_n(x) &= \int_{\mathcal X} \kappa(x,x') \, d\mu_{X,n}(x'),\\
    v_n(x) &= \int_{\mathcal X}\frac{\kappa(x,x')}{u_n^\alpha(x')} \, d\mu_{X,n}(x'). 
\end{align*}
Here, $\kappa$ is a continuous, symmetric, strictly positive, positive-definite kernel (e.g., the radial Gaussian kernel from~\eqref{eqKGauss}), $u_n$ is the normalization function from~\eqref{eqUV}, and $\alpha $ is a real parameter (typically set to 0, $1/2$, or 1). One can verify that the kernel $p_n$ just defined satisfies the detailed balance condition in~\eqref{eqBalance} for $w=p_n$ and $ d =  v_n / u_n^\alpha $, and thus can be employed for KAF as described in Section~\ref{secExt}. An advantage of these kernels over the kernels in~\eqref{eqBistoch} is that they do not require integration with respect to $x''$ in their definition, thus avoiding a source of sampling error. A disadvantage is that one needs to keep track of a biorthogonal pair of bases, $\{ \xi_i \}$ and $\{ \xi'_i \}$, as opposed to a single orthonormal basis $\{ \phi_i \}$ in the symmetric case. 
    
\subsection{\label{secDelays}Kernels based on delay-coordinate maps}

Delay-coordinate maps is a technique originally introduced for empirical state space reconstruction of partially observed dynamical systems \cite{PackardEtAl80,Takens81}, which employs the time ordering of the covariate data to embed it in a higher-dimensional space. Choosing a positive integer parameter $q$ (the number of delays), we define the covariate map $X_q : \Omega \to \mathcal X^q$ such that
\begin{equation}
    \label{eqDelay}
    X_q(\omega) = (X(\omega), X(\Phi^{-\Delta t}(\omega)), \ldots, X(\Phi^{-(q-1)\,\Delta t}(\omega))).
\end{equation}
Note that $X_q(\omega)$ can be empirically evaluated without explicit knowledge of the dynamical flow $\Phi^t$ on state space, so-long as time-ordered covariate data $X(\Phi^{-j,\Delta t}(\omega))$ are available over the temporal window $[-(q-1)\,\Delta t, 0 ]$ at the sampling interval $\Delta t$. 

Intuitively, this approach should increase the information content of the covariate map, since the time ordering of the data is a manifestation of the underlying dynamical flow. This intuition has been made mathematically precise in a number of ``embedology'' results for flows in finite-dimensional state spaces \cite{SauerEtAl91}, as well as classes of partial differential equation models \cite{Robinson05}. These studies have established that under natural assumptions, $X_q$ becomes an injective map for sufficiently large $q$, even if the raw covariate map $X$ is not injective. In that case, the temporal evolution of the covariate $X_q(\omega)$ on the support $ \mu \circ X_q^{-1}$ of the invariant measure on delay-coordinate space becomes a homeomorphic copy of the dynamics on the support $\Omega_\mu$, with optimal potential predictability. While rigorously verifying the appropriate embedding conditions in practice is oftentimes difficult, it is generally expected that including delays can recover at least some of the state space degrees of freedom lost due to non-injectivity of $X$, leading to more skillful forecasts. Indeed, delay-coordinate maps have been employed in a number of parametric \cite{SmallJudd98} and nonparametric  \cite{FarmerSidorovich87,SugiharaMay90,Sauer92} forecasting methodologies, and have also been found useful for extraction of coherent features from time series data \cite{BroomheadKing86,VautardGhil89,GiannakisMajda12a,BerryEtAl13,ArbabiMezic17,BruntonEtAl17,DasGiannakis19,Giannakis19}.

In the context of KAF, we generally expect the conditional expectation $\mathbb E[U^\tau Y \mid X_q]$ approximated by the algorithm to exhibit smaller intrinsic error (see~\eqref{eqSigmaRho}) with increasing $q$, and thus smaller total error for appropriately constructed hypothesis spaces. One way of constructing these spaces is to employ the strictly positive-definite, Markovian kernels described in Sections~\ref{secPD} and~\ref{secMarkov}, replacing the ``snapshot'' data in $\mathcal X$ with delay-embedded sequences in $\mathcal X^q$. For instance, a natural analog of the local kernel in~\eqref{eqKLocal} is $k^{(q)} : \mathcal X^q \times \mathcal X^q \to \mathbb R  $, where
\begin{equation}
    \label{eqKLocalQ}
    k^{(q)}(\tilde x, \tilde x' ) = h( L_q(\tilde x,\tilde x') / \epsilon).
\end{equation}
Here, $L_q$ is the symmetric function 
\begin{displaymath}
    L_q(\tilde x, \tilde x') = \frac{1}{q} \sum_{j=0}^{q-1} L(x_j,x'_j), 
\end{displaymath}
induced on $\mathcal X^q$ from $L$, where $ \tilde x = ( x_1, \ldots, x_q ) $ and $\tilde x' = ( x'_1, \ldots, x'_q )$. The kernel $k^{(q)}$ can then be normalized via the symmetric or non-symmetric normalization procedures described in Section~\ref{secMarkov} to yield a Markov kernel.  See, e.g., \cite{ComeauEtAl17,ComeauEtAl19,AlexanderEtAl17,WangEtAl19b} for applications of KAF with Gaussian kernels on delay-coordinate space.

Before closing this section, we should point out that while beneficial from the point of view of topological state-space reconstruction, delay-coordinate maps with large numbers of delays may face potential limitation from the point of view of spectral characteristics of the underlying dynamical system. In particular, observe that the pullback $L_{\Omega,q}: \Omega \times \Omega \to \mathbb R_+$ of $L_q$ on state space, i.e., $L_{\Omega,q}(\omega,\omega') = L_q(X_q(\omega),X_q(\omega')) $ has the structure of an ergodic average of the continuous function $L_\Omega(\omega,\omega') = L(X(\omega),X(\omega'))$ under the product dynamical flow $ \Phi^t \times \Phi^t$, viz.
\begin{displaymath}
    L_{\Omega,q}(\omega,\omega') = \frac{1}{q} \sum_{j=0}^{q-1} L_\Omega(\omega_{-j}, \omega'_{-j} ), 
\end{displaymath}
where $\omega_j = \Phi^{j\,\Delta t}(\omega)$ and $\omega'_j = \Phi^{j\,\Delta t}(\omega')$. As a result, by the Birkhoff pointwise ergodic theorem \cite{Walters81}, as $q\to \infty$, $L_{\Omega,q}$ converges $\mu \times \mu$-almost everywhere to a function $L_{\Omega,\infty} \in \mathbb L^2(\mu \times \mu)$. Further, it can be shown \cite{DasGiannakis19} that in the same limit, the kernel integral operators $G_{\Omega,q}$ on $L^2(\mu)$ associated with the pullback kernel $k^{(q)}_\Omega : \Omega \times \Omega \to \mathbb R$ induced by~\eqref{eqKLocalQ}, $k_{\Omega}^{(q)}(\omega,\omega')=k^{(q)}(X(\omega),X(
\omega'))$ converge in operator norm, and thus in spectrum, to a Hilbert-Schmidt integral operator $ G_{\Omega,\infty} : L^2(\mu) \to L^2(\mu) $ associated with the kernel 
\begin{displaymath}
    k_{\Omega,\infty}(\omega,\omega') = h(L_{\Omega,\infty}(\omega,\omega')/\epsilon).
\end{displaymath}

Now, by invariance of Birkhoff averages, at any time $ t \in \mathbb R$ the kernel $k_{\Omega,\infty}$ is invariant under the Koopman operator $U^t \otimes U^t : L^2(\mu \times \mu) \to L^2(\mu\times \mu) $ for the product system, and the latter implies that $ U^\tau $ and $ G_{\Omega,\infty}$ are commuting operators \cite{DasGiannakis19}. As a result, every eigenspace of $G_{\Omega,\infty} $ at nonzero eigenvalue  (which is finite-dimensional by compactness of this operator) is a finite union of Koopman eigenspaces. The latter implies, in particular, that the nullspace of $G_{\Omega,\infty} $ must necessarily contain the Koopman-invariant subspace  $ L^2(\mu)$ associated with the continuous spectrum of $U^\tau $ (see \cite{Halmos56,Mezic05,DasGiannakis19} for precise definitions of this subspace). If it now happens that the subspace $L^2_X(\mu)$ where the conditional expectation $ \mathbb E[ U^\tau Y \mid X ] $ lies has a nonzero intersection with the continuous-spectrum subspace, as would typically be the case in systems with sufficiently complex (mixing) dynamics, then the hypothesis spaces associated with $k_{\Omega,\infty} $ will fail to be dense in $L^2_X(\mu)$, and thus the KAF target functions may fail to converge to $\mathbb E[U^\tau Y \mid X]$. Since the empirical integral operators with large numbers of delays are spectrally close to $ G_{\Omega,\infty}$, this behavior indicates that there may be situations where increasing $q$ beyond a limit may be detrimental to forecast skill.

\section{\label{secApps}Applications}
We present two examples to illustrate how to build a kernel forecasting function, as well as some basic properties of convergence to the conditional expectation. See Table~\ref{tablePseudocode} for a summary of the algorithmic steps involved in KAF.

\subsection{\label{secCircle}Circle rotation}

Our first example is periodic flow on the circle, $\Omega = S^1$. Expressed in terms of canonical angle coordinates $ \omega \in [ 0, 2 \pi )$, the dynamical flow map $\Phi^t$ takes the form of a translation, 
\begin{displaymath}
    \Phi^t( \omega ) = \omega + \alpha t \mod 2 \pi, \quad \alpha \in \mathbb{R},
\end{displaymath}
with a period of $ 2 \pi/\alpha $, exhibiting a unique ergodic invariant Borel probability measure $ \mu $, equal to a normalized Lebesgue measure. As covariate and response spaces, we choose $ \mathcal{X} = \mathcal{Y} = \mathbb{R}$, and we prescribe covariate and response maps, $X$ and $Y$, respectively, given by simple trigonometric functions as follows: 
\begin{displaymath}
    X(\omega) = \cos(\omega), \quad Y(\omega) = \sin(\omega).
\end{displaymath}

Under this setup, we have $U^\tau Y(\omega) = \sin( \omega+\alpha\tau)$, and the conditional expectation $\mathbb{E}[U^\tau Y \mid X = x]$ is the average of $U^\tau Y$ at the two angles for which $X(\omega+ \alpha\tau) =  x$; specifically,
\begin{align*}
    Z_\tau(x) &= \frac{ \sin( \arccos(x) + \alpha\tau) + \sin(- \arccos(x) + \alpha\tau)}{2} \\
    &= x\sin(\alpha \tau).
\end{align*}
The intrinsic error $\sigma_\tau$ may then be computed directly as
\begin{displaymath}
    \sigma_\tau = \lVert Y \rVert_{L^2(\mu)}^2 \cos^2(\alpha\tau), \quad \lVert Y \rVert_{L^2(\mu)} = 1 / \sqrt{2}.
\end{displaymath}
Observe that the intrinsic error is maximal (and equal to the squared $L^2(\mu)$ norm of the response variable) when $\tau = q \pi / \alpha$, and minimal (and equal to zero)  when $\alpha\tau = (2q+1)\pi/2$, where $ q $ is any integer.

The pushforward measure $ \mu_X $ in covariate space is supported on the interval $ [ -1, 1 ] \subset \mathbb{R} $, where it has the density 
\begin{displaymath}
    \varrho(x) = \frac{d \mu_X(x)}{d\Leb} = \frac{1}{2\pi \sqrt{1-x^2}}
\end{displaymath}
relative to Lebesgue measure. Note that $ \varrho(x) $ diverges at the boundary points $ x = \pm 1 $, but nevertheless lies in the $\mathbb{L}^1$ space associated with the Lebesgue measure. Given a kernel $ k : \mathbb{R} \times \mathbb{R} \to \mathbb{R} $ meeting the conditions of Section~\ref{secError}, the eigenvalue problem for the associated integral operator $ G $ then becomes 
\begin{displaymath}
\int_{-1}^1 k(x,x') \phi_i(x') \varrho(x') \, dx' = \lambda_i \phi_i(x).
\end{displaymath}

A closed, analytic expression for this eigenvalue problem is not known for arbitrary choices of kernel $k$. Instead, using the radial Gaussian kernel from~\eqref{eqKGauss}, we numerically solve the eigenvalue problem for the data-driven operator $G_n $, constructed from a sequence $ x_1, \ldots, x_n \in \mathbb{R} $ of covariate points obtained from an underlying dynamical trajectory $ \omega_1, \ldots, \omega_n \in S^1 $, as described in Section~\ref{secTarget}. Using the corresponding sequence $ y_{1}, \ldots, y_{n} \in \mathbb{R} $ of response variables, we then build the KAF target function $f_{\tau,\ell,n}$ via~\eqref{eqKAF}.

Here, we set the frequency $\alpha = \sqrt{2}$, and employ a training dataset of $ n = 1000$ samples, taken at an interval of $ \Delta t = 2\pi/100 $ time units. Note that $\Delta t $ is rationally independent from the rotation period, which ensures that the discrete-time map $ \Phi^{\Delta t} $ provides an ergodic sampling of $ \mu $. Using this dataset, we have computed $f_{\tau,\ell,n} $ for lead times $ \tau = q \, \Delta t / \alpha $, with $ q $ an integer in the interval $ [ 0, 24 ] $. To assess forecast skill, one can compute the mean square error (MSE)
\begin{displaymath}
    \tilde{\mathcal{E}}_{\tau,m}(f_{\tau,\ell,n}) = \frac{1}{m} \sum_{j=1}^m \lvert f_{\tau,\ell,n}( \tilde x_j ) - U^\tau Y( \tilde \omega_j ) \rvert^2
\end{displaymath}
on a verification dataset $ \tilde \omega_j = \Phi^{(j-1)\, \Delta t}(\tilde \omega_1) $, $ \tilde x_j = X(\tilde \omega_j ) $. Since the intrinsic error happens to be analytically expressible for this problem, we report the empirical excess generalization error
\begin{displaymath}
    \tilde{\mA}_{\tau,m}(f_{\tau,\ell,n}) = \frac{1}{m} \sum_{j=1}^m \lvert f_{\tau, \ell, n}( \tilde x_j ) - Z_\tau( \tilde x_j ) \rvert^2. 
\end{displaymath}
The empirical MSE and excess generalization error approximate the true MSE and excess generalization error, $ \mathcal{E}_\tau(f_{\tau,\ell,n})$ and $\mathcal{A}_\tau(f_{\tau,\ell,n})$, respectively, and converge to these quantities as $m \to \infty$. Here, we employ a verification dataset of  $m = \text{10,000}$ samples, starting from a state $ \tilde \omega_1 $ chosen randomly and uniformly on the interval $[0,2\pi)$ (so that $ \omega_1 - \tilde \omega_1$ and $ \Delta t$ are rationally independent with probability 1).   
\begin{figure*}
    \centering\includegraphics[width=.7\linewidth]{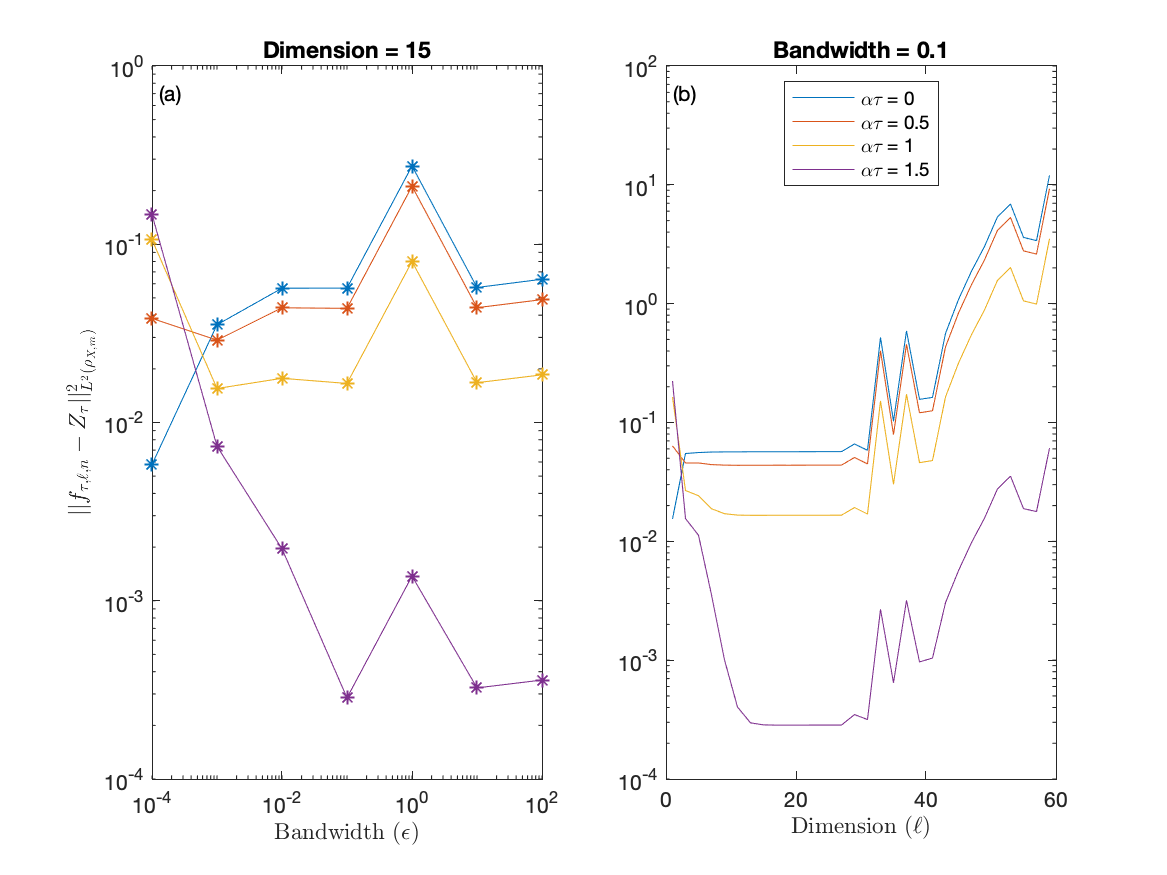}
\caption{\label{figCircle} Empirical excess generalization error $\tilde{\mathcal{A}}_{\tau,m}(f_{\tau,\ell,n})$ of KAF target functions $f_{\tau,\ell,n}$ for periodic flow on the circle with angular frequency $ \alpha = \sqrt{2} $, obtained using  the radial Gaussian kernel and a training dataset of $n = 1000 $ samples. The empirical excess generalization error is plotted for the lead times $\tau \in \left\{0, 1/2\alpha, 1/\alpha, 3/2\alpha\right\}$, using $m = \text{10,000}$ test points. In Panel~(a), the dimension of the hypothesis space is fixed at $\ell = 15$, and the bandwidths are varied to take on values $\epsilon \in \left\{ 10^{-4}, 10^{-3}, \dots, 10^2 \right\} $. In Panel~(b), the bandwidth is fixed at $\epsilon = 0.1$, and the dimensions take on all odd integers from $1$ to $59$. For reference, note that the squared $L^2(\mu)$ norm of the response variable $U^\tau Y $  is equal to 1/2.}
\end{figure*}

Figure~\ref{figCircle} shows the absolute value of empirical excess generalization error, plotted against the bandwidth parameter $\epsilon$ of the Gaussian kernel and the hypothesis space dimension $\ell$ for representative lead times $ \tau $ in the range $ [ 0, \pi / ( 2 \alpha ) ] $. In Figure~\ref{figCircle}(a), $\ell$ is kept fixed at 15, and $\epsilon $ varies logarithmically in the interval $ [ 10^{-4}, 10^2 ] $. The results show agreement between several different choices of $ \epsilon$ in some regimes of $\tau$, but also notable discrepancy in the region where the intrinsic error is already very small (i.e., when $\alpha\tau$ is close to $\pi/2$). In such a regime, the less sensitive kernels of large bandwidth are better able to capture that the generalization error is close to 0. In general, the $\tilde{\mathcal{A}}_{\tau,m}$ values in Figure~\ref{figCircle}(a) lie approximately in the interval $ [ 10^{-4}, 3 \times 10^{-1} ] $, which corresponds to approximately $ 2 \times 10^{-4} $ to $ 5 \times 10^{-1} $ multiples of the squared $L^2(\mu) $ norm of the covariate variable.

Figure~\ref{figCircle}(b) shows the behavior of empirical excess generalization error at fixed $ \epsilon = 0.1 $ and representative values of $ \ell $ in the range 1 to 60. Employing just the first eigenfunction performs best for $\alpha \tau = 0$, but employing more eigenfunctions is better for larger values of $\alpha \tau$.  Most notable, however, is the characteristic bias-variance tradeoff as $\ell$ increases, with a valley of optimal values of $\ell$ between 10 and 30. For instance, at $\alpha \tau = 1.5 $, the error decreases from $ \simeq 2 \times 10^{-1} $ for $\ell = 1 $ to a minimal value of $ \simeq 2 \times 10^{-4} $ for $ \ell = 20 $, but then increases for larger $ \ell $ to $ \simeq 10^{-2} $ values. This is a manifestation of the fact that the true error $\mathcal{E}_\tau(f_{\tau,\ell,n})$ may increase with $ \ell $ at fixed $ \tau $ and $ n $, even though the empirical error $\mathcal{E}_{\tau,n}(f_{\tau,\ell,n})$ is always a non-increasing function of $ \ell $.

\subsection{\label{secL63}Lorenz 63 system}

In the L63 system \cite{Lorenz63}, the state space is $\Omega = \BR^3$. The dynamical flow $\Phi^t( \omega_0) $ starting from $ \omega_0 \in \mathbb{R}^3$ is given by solution of the initial-value problem 
\begin{displaymath}
    \dot \omega(t) = \vec V(\omega(t)), \quad \omega(0) = \omega_0, 
\end{displaymath}
where $ \vec V : \mathbb{R}^3 \to \mathbb{R}^3 $ is the smooth vector field with components $ ( V^1, V^2, V^3 )  $ at $ \omega = (\omega^1,\omega^2,\omega^3) $ given by $V^1 = \sigma( \omega^2-\omega^1) $, $ V^2=\omega^1(\mu-\omega^3) $, and $ V^3 = \omega^1 \omega^2 - \beta \omega^3$. Here, $\beta$, $\mu$, and $\sigma$ are real parameters, which we set to the classical values $ \beta = 8/3 $, $\mu = 28$, and $\sigma = 10 $. For this choice of parameters, the L63 system is rigorously known to have a compact attractor $ \Omega_\mu \subset \mathbb{R}^3  $  \cite{Tucker99} with fractal dimension $\approx 2.06$ \cite{McGuinness68}, supporting a physical invariant measure $ \mu$ with a single positive Lyapunov exponent $ \Lambda \approx 0.91$ \cite{Sprott03}. Due to dissipative dynamics, the attractor is contained within absorbing balls \cite{LawEtAl14}, ensuring the existence of the compact set  $\mathcal{U} \subseteq \mathcal{X}$ in covariate space. In light of these facts, all of the assumptions on the dynamical system  made in Section~\ref{secSampleError} rigorously hold. The L63 system is also rigorously known to be mixing \cite{LuzzattoEtAl05}, and thus exhibits the loss of long-term predictability discussed in Section~\ref{secMixing}.

In the experiments that follow, we shall let $\mY = \BR$, and let the response variable $Y : \Omega \to \mY$ pick out one of the state vector components, i.e., if $\omega = (\omega^1,\omega^2,\omega^3) \in \mathbb{R}^3$, then $Y(\omega) = \omega^i$, for either $i = 1, 2, 3$. To illustrate the conditional probability framework discussed in Section~\ref{secProb}, for each such response variable we will consider the event $ S = \{ \omega \in \Omega : Y(\omega) >\theta \} $, where $ \theta $ is an empirical mean of $ Y $ computed from the training data. That is, we will use KAF to estimate the conditional probability that the components of the state vector exceed their mean values.  As for the covariate variable $X$, we will consider two cases, namely,  full observations, $\mathcal{X} = \mathbb{R}^3$ and $ X = \Id $, and a partially observed setup with $ \mathcal{X} = \mathbb{R} $ and $X(\omega) = \omega^1 $.  Consequently, in the partially observed setup  $\BE[U^\tau Y \mid X]$ represents the conditional expectation of the $i$-th component of $\Phi^\tau(\omega) $, given the first component of $ \omega $. 

All experiments use covariate data $ x_j = X(\omega_j) $ and response data $ y_{\tau,j} = U^\tau Y( \omega_j ) $ generated from the same underlying trajectory $ \omega_1, \ldots, \omega_n \in \mathbb{R}^3$, with $ \omega_j = \Phi^{(j-1) \, \Delta t}( \omega_j ) $.  The trajectory $ \omega_j $ was numerically generated in Matlab using the \texttt{ode45} solver, starting from an arbitrary initial condition $\omega_0$ and waiting for a long spinup time before collecting the first sample $ \omega_1 $. We nominally work with a training dataset consisting of $ n =\text{64,000} $ samples, taken at a sampling interval $ \Delta t $ equal to 0.01 natural time units (i.e., about $1/100 $ of the characteristic Lyapunov timescale $ 1/ \Lambda $ of the system). Additional experiments with dataset sizes ranging from $n= 640$ to $n = \text{512,000} $ and/or a longer sampling interval of $\Delta t = 0.1$ were also conducted to investigate the performance of KAF relative to the $(n,\Delta t) = (\text{64,000},0.01)$ baseline case. 
    
To assess forecast skill, we use empirical root mean square error (RMSE) metrics computed from an independent verification dataset as in the circle example of Section~\ref{secCircle}. Specifically, the RMSE of the target function $f_{\tau,\ell,n}$ at lead time $ \tau $ is given by $ \sqrt{\tilde{\mathcal{E}}_{\tau,m}(f_{\tau,\ell,n})}$, where the verification dataset has the same number of $m = \text{64,000} $  samples as the training dataset, and was  obtained via a similar spinup procedure starting from a different initial condition. Note that aside from the covariate and response data in the training phase, and the covariate data in the verification phase, no other information about the system state and/or dynamics was provided to the KAF algorithm.

\begin{figure*}
    \includegraphics[width=\linewidth]{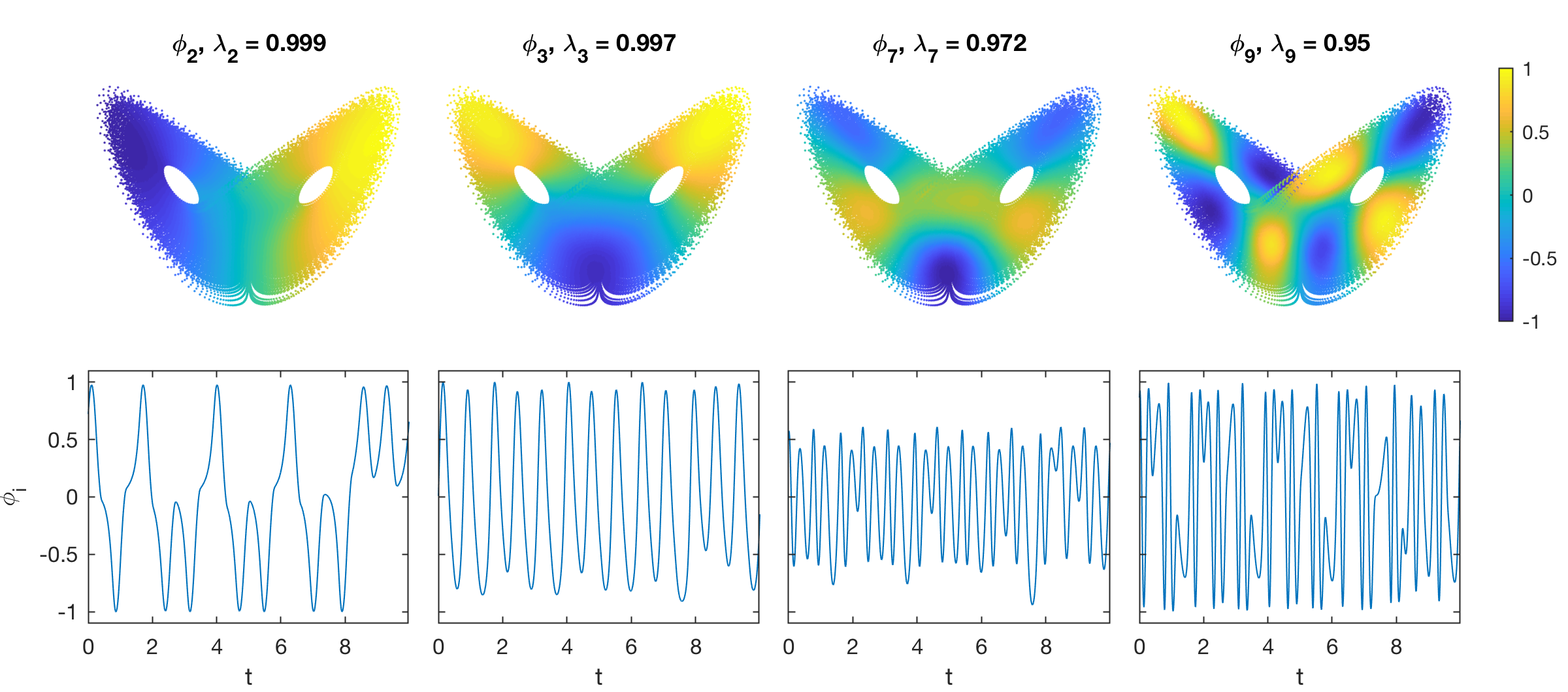}
    \caption{\label{figL63PhiFull}Representative data-driven eigenfunctions $\phi_{i,n}$ and their corresponding eigenvalues $\lambda_{i,n}$, computed from the fully observed L63 system. Top: Scatterplots of the eigenfunction values $ \phi_{i,n}(x_j) $ on the covariate training data $x_j = \omega_j \in \mathbb{R}^3$. Bottom: Eigenfunction time series $t_j \mapsto \phi_{i,n}(x_j)$ over a portion of the training dataset spanning 10 natural time units. Notice that, despite the fact that the L63 attractor is not a Riemannian manifold, the eigenfunctions qualitatively resemble a generalized Fourier basis associated with a heat kernel. That is, as $\lambda_{i,n}$ decreases, $\phi_{i,n}$ exhibits increasingly small-scale oscillatory behavior, allowing one to represent functions of increasingly fine structure through eigenfunction expansions.}  
\end{figure*}

\begin{figure*}
    \includegraphics[width=\linewidth]{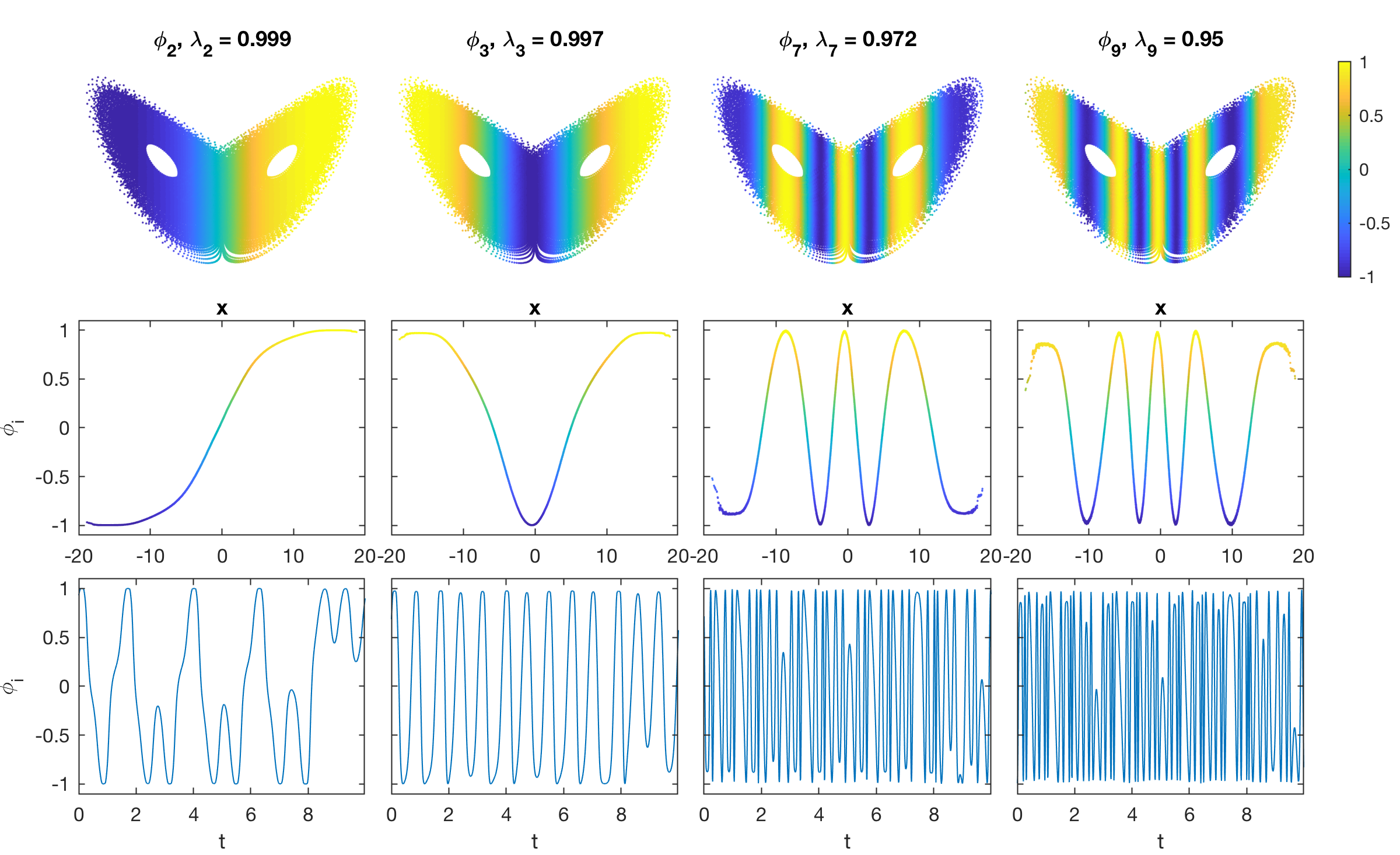}
    \caption{\label{figL63PhiPartial}As in Figure~\ref{figL63PhiFull}, but for eigenfunctions computed from the partially observed L63 system. Top: Scatterplots of the eigenfunction values $ \phi_{i,n}(x_j) $ on the L63 states $ \omega_j $ underlying the covariate data $x_j \in \mathbb{R}$. Middle: Eigenfunction values $\phi_{i,n}(x_j)$ versus $x_j$. Coloring is as in the top row. Bottom: Eigenfunction time series $t_j \mapsto \phi_{i,n}(x_j)$ over a portion of the training dataset spanning 10 natural time units.}
\end{figure*}

The first step in the KAF pipeline  is to compute the kernel eigenfunctions $ \phi_{i,n} $, whose corresponding RKHS functions $ \psi_{i,n} $ form orthonormal bases for the hypothesis spaces $\mathcal{H}_{n,\ell}$. For that, we employ the variable-bandwidth, Markov-normalized kernels from~\eqref{eqBistoch}, with automatically tuned bandwidth and dimension parameters (see Section~\ref{secKernel}). Representative eigenfunctions for the fully- and partially-observed systems are displayed in Figures~\ref{figL63PhiFull} and~\ref{figL63PhiPartial}, respectively, in both scatterplot and time series form. There, it is evident that the eigenfunctions behave like a generalized Fourier basis on the support of the measure $ \mu_X$, with eigenfunctions corresponding to smaller eigenvalues allowing to resolve functions of increasingly smaller-scale variability on the L63 attractor. Notice, in particular, that in the partially observed example with $\mathcal{X} = \mathbb{R}$, the $ \phi_{i,n} $ are increasingly oscillatory, orthogonal functions on the real line, which pull back to $\mathcal{G}$-measurable functions on the attractor in $\mathbb{R}^3$ with no variability in the $\omega^2 $ and $ \omega^3 $ coordinates. It is precisely such a lack of variability that contributes to degraded forecast skill when faced with non-injective covariate functions.    

Next, using the eigenfunctions and the response data in the training phase, we construct the empirical target functions $f_{\tau,\ell,n}$ from~\eqref{eqKAF}. We also compute the error estimators $\varepsilon_{\tau,\ell,n}$ from~\eqref{eqErrEst}, which we use to place ``error bars'' around our forecasts of the form $ f_{\tau,\ell,n}(x) \pm \varepsilon_{\tau,\ell,n}(x)$. According to Section~\ref{secProb}, for an unbiased error estimate, the RMS value of $\varepsilon_{\tau,\ell,n}$ in the verification phase should be close to the actual RMS forecast error. We use $\ell = 3000 $ and 1000 eigenfunctions for the fully observed and partially observed setups, respectively.

\paragraph{Forecasting state vector components} Figure~\ref{figL63Pred} shows prediction results for the three components of the state vector for the fully observed and partially observed systems, together with error estimates based on $ \varepsilon_{\tau,\ell,n} $. We show representative forecast trajectories starting from an arbitrary initial condition in the verification dataset, as well as aggregate RMSE scores as a function of lead time, normalized by empirical standard deviation (i.e., the $L^2(\mu_n)$ norm of $Y - \int_\Omega Y \, d\mu_n$). Starting from the fully observed examples,  the RMSE of all three state vector components $ \omega^i $ exhibits an initial exponential-like increase from near-zero values for $\tau \lesssim 0.5 \simeq 0.5 / \Lambda $. This period is followed by an intermediate-time regime with  more gradual RMSE increase and noticeable oscillatory behavior,  until convergence to the equilibrium standard deviation (normalized RMSE $\simeq 1$) at late times, $\tau \gtrsim 4$. 

Examining the individual forecast trajectories, it is evident that the late-time convergence of the RMSE to a near-constant values is a manifestation of the trajectories converging to the mean, $\mathbb{E}[Y]$. The numerical results are therefore consistent with the theoretically expected late-time behavior of KAF in the presence of mixing dynamics, discussed in Section~\ref{secMixing}. It is also evident from Figure~\ref{figL63Pred} that the error estimators $ \varepsilon_{\tau,\ell,n} $ provide useful uncertainty quantification. That is, the error bars derived from these quantities envelop, for the most part, the true trajectories, and their RMS values agree well with the forecast RMSE. 

Overall, in the fully-observed experiments, $\omega^3$ is the most predictable state vector component (likely due to symmetry of the L63 equations), followed by $ \omega^1 $ and $ \omega^2 $ which are nearly equally predictable (again due to symmetry). If one were to set a normalized RMSE value of 0.6 as a threshold for loss of skill, $ \omega^3 $ would remain predictable out to $\simeq 3 $ natural time units (i.e., $\simeq 3$ Lyapunov timescales), whereas $ \omega^1 $ and $ \omega^2 $ would remain predictable out to $\tau \simeq 2 $. Setting that threshold to 0.8 increases the predictability horizon of $ \omega^3 $ and $ \omega^1$/$\omega^2$ to $ \tau \simeq 5 $ and 2.75, respectively.  

Turning now to the Figure~\ref{figL63Pred} results for the partially observed system, it is clear that the act of observing $\omega^1$ only in the covariate space bears a significant impact on forecast skill, particularly for $ \omega^2 $ and $ \omega^3$. Indeed, for these state vector components, the non-injectivity of the covariate function means that the normalized RMSE can be significant even at $ \tau = 0 $, without ever dropping below $ \lesssim 0.4 $. Yet, even though the method cannot overcome the intrinsic error of this observational setup, it is nevertheless capable of providing fairly adequate uncertainty quantification, as manifested by the reasonably good ability of the estimated error bars to envelop the true trajectories (with the notable exception of certain extremal points) and the close agreement between the RMS values of $\varepsilon_{\tau,\ell,n}$ and the forecast RMSE.          

\begin{figure*}
    \includegraphics[width=\linewidth]{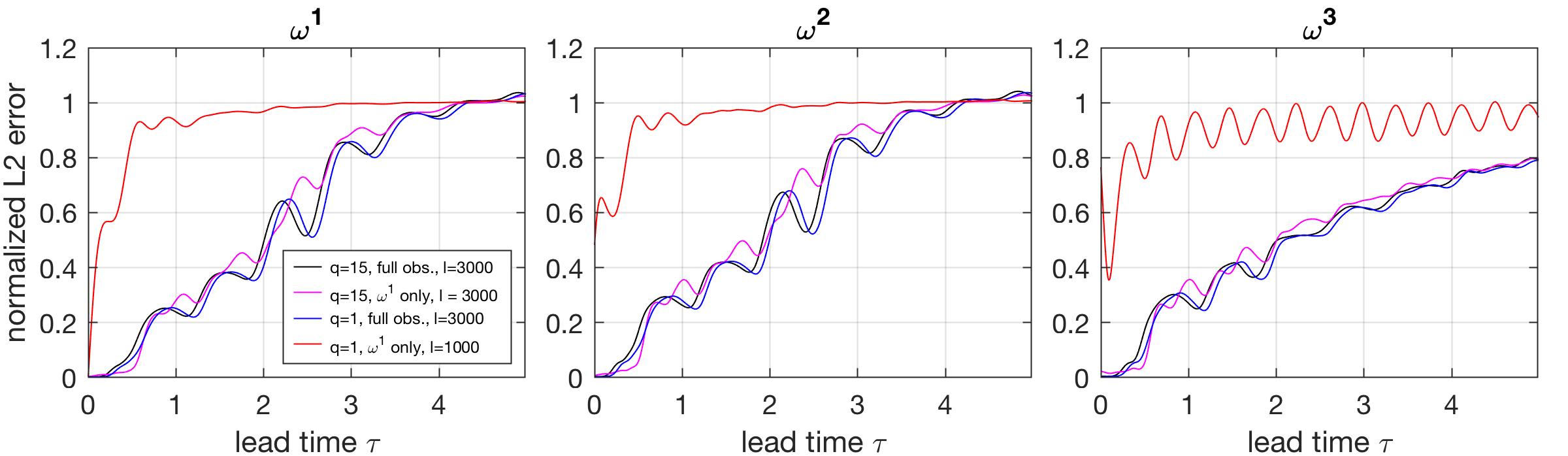}
    \caption{\label{figL63RMSE_delays}Normalized RMSE versus lead time plots for KAF applied to the state vector components of the L63 system, illustrating the effects of incorporating delays in the covariate function. The black and red lines show the RMSE scores for the fully and partially observed cases from Figure~\ref{figL63Pred}, respectively, using $n = 64,000 $ training samples and $\ell = 3000$ (full observations) and $\ell = 1000$ (partial observations) eigenfunctions. The blue and magenta lines show RMSE scores obtaining by incorporating $ q= 15$ delays to the full and partial observation maps, respectively. The generally high consistency of the $q=15 $ results with those of the fully observed experiments without delays illustrate that (i) delay-coordinate maps successfully recover information loss due to partial observations; and (ii) KAF behaves stably in delay-spaces with potentially poor conditioning due to colinearity of delay coordinates.}
\end{figure*}

\paragraph{Delay-coordinate maps} As a demonstration of the efficacy of delay-coordinate maps in recovering forecast skill lost due to forecast observations, in Figure~\ref{figL63RMSE_delays} we compare the RMSE scores from the fully and partially observed experiments in Figure~\ref{figL63Pred} with their counterparts obtained by including $q=15$ delays to the respective covariate maps. Specifically, we construct delay-coordinate maps via~\eqref{eqDelay} based on either $X(\omega)=(\omega^1,\omega^2,\omega^3)$ or $X(\omega) = \omega^1$, and build KAF models using the same class of kernels, sampling intervals, and training dataset sizes as the experiments without delays. In both cases we use the same number of eigenfunctions as the fully observed case with no delays, $\ell = 3000$. Note that the delay-coordinate experiments based on full observations are interesting despite the fact that there is no potential predictability to be gained---this is because incorporating delays can introduce strong colinearities in the training data, increasing the likelihood of overfits. In addition, the L63 system falls squarely in the class of mixing dynamical systems discussed in Section~\ref{secDelays}, where incorporating delays in the kernel can suppress the nonzero eigenvalues of the corresponding integral operators, making them prone to sampling errors. 

As is evident in Figure~\ref{figL63RMSE_delays}, adding delays results in a considerable increase of  skill for the partially observed system, with RMSE scores in generally good agreement with the fully observed systems. Interestingly, the RMSE values for short-term forecasts with $ 0.2 \lesssim \tau \lesssim 0.5$ appear to be smaller than the fully observed experiments (either with or without delays), though longer-term forecasts exhibit intervals (e.g., $ \tau \in [2.2, 2.5]$) where the RMSE of the partially observed system with delays is noticeably higher. Meanwhile, the fully observed system with delays exhibits very comparable skill as the system without delays, despite the issues associated with colinearity and mixing dynamics mentioned above.    

\begin{figure*}
    \includegraphics[width=\linewidth]{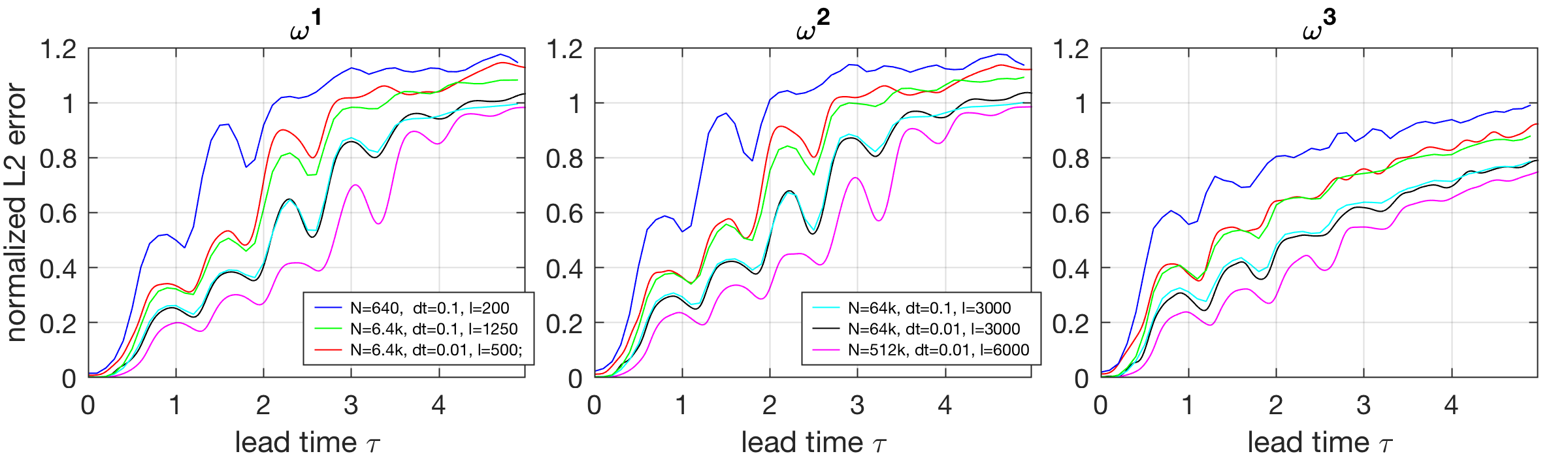}
    \caption{\label{figL63RMSE}Normalized RMSE versus lead time plots for KAF applied to the state vector components of the L63 system, under full observations and for a variety of numbers of training samples ($n$), sampling intervals ($\Delta t$), and hypothesis space dimensions ($\ell$). The case with $n = \text{64,000}$ and $\Delta t = 0.01$ shown here in a black line is identical to the fully observed case from Figure~\ref{figL63Pred}. The other experiments shown here in colored lines correspond to various modifications of sample number and/or sampling frequency relative to the $ n = \text{64,000}$, $\Delta t = 0.01$ baseline case.}
\end{figure*}

\paragraph{Sensitivity analysis} The results displayed in Figure~\ref{figL63Pred} were obtained using a fairly dense sampling of the L63 attractor, corresponding to $n \, \Delta t = 640$ natural time units, or, approximately 800 oscillations assuming a characteristic oscillatory timescale of 0.8. Moreover, the sampling interval $\Delta t = 0.01$ was short compared to the oscillation and Lyapunov timescale of the system. To assess the performance of KAF in environments with shorter and less frequently sampled training data, we have performed a suite of forecasting experiments with full observations ($X=\Id$) that differ from our nominal setup with $n = \text{64,000}$ and $\Delta t = 0.01$ by various modifications of the number of training samples $n$ and sampling interval $ \Delta t $. RMSE results from these experiments are depicted in Figure~\ref{figL63RMSE}.

Before delving to a description of these results, let us recall that, as with any supervised learning technique, KAF strives for a balance between generalization error (the difference between the ideal target function $f_{\tau,\ell}$ and the regression function) and sample error (the difference between $f_{\tau,\ell}$ and the empirical target function $ f_{\tau,\ell,n}$). This balance is attained by controlling the number of eigenfunctions (principal components) $\ell$ employed, and the generalization error is a decreasing function of $\ell$. On the other hand, at fixed $n$, the sample error is generally an increasing function of $\ell$. As the sampling provided by the training data becomes poorer (by decreasing $n$ and/or $\Delta t$), the values of $\ell$ achieving that balance become smaller, generally resulting to a decrease of forecast skill. 

The RMSE results in Figure~\ref{figL63RMSE} are for values of $\ell$ chosen on the basis of yielding good skill over the full range $\tau \in [ 0, 5 ] $ of lead times examined. It should be noted that we did not perform an exhaustive search to select these values, as we found that the dependence of skill on $\ell$ exhibits plateau behavior, analogously to the circle example in Figure~\ref{figCircle}(b). We should also point out that in an ``operational'' environment one would typically select different values of $\ell$ for each lead time $\tau$ so as to minimize RMSE. In particular, that for a mixing dynamical system such as L63, as $\tau $ increases the conditional expectation $\mathbb E[U^\tau Y \mid X]$ weak-converges to a constant, indicating that for the class of Markov kernels employed in this work (where the top eigenspace is spanned by constant functions) smaller $\ell$ values may be warranted at large $\tau$.

With these considerations in mind, we now turn to the results in Figure~\ref{figL63RMSE}. First, note that increasing the sampling interval by a factor of 10 to $\Delta t = 0.1$, while keeping the number of training samples and eigenfunctions fixed to the nominal values from Figure~\ref{figL63Pred}, $(n,\ell) = (\text{64,000},3000)$, imparts little change to forecast skill. This suggests that the convergence of the leading $\ell$ eigenspaces of the empirical integral operators $G_n $ employed in KAF is largely unaffected by this reduction of sampling frequency. On the other hand, a tenfold reduction of the number of samples to $n = 6400$, using either $\Delta t = 0.01$ or $\Delta t = 0.1$, does impart a noticeable reduction of skill, as we are forced to work with smaller numbers of eigenfunctions, $\ell = 500  $ and $1250$ for $ \Delta t = 0.01$ and 0.1, respectively. Nevertheless, at least over short to moderate lead times, $\tau \lesssim 2$, the reduction of skill compared to the $n = \text{64,000}$ cases is fairly modest. For instance, using again a 0.6 value of normalized RMSE as a useful-skill threshold, the forecasts of the $\omega^1$ and $\omega^2$ variables based on the  $n = 6400$ remain useful out to $\tau \simeq 1.8$ versus $ \tau \simeq 2.1$ for $n = \text{64,000}$. Larger discrepancies are observed at longer leads, $\tau \gtrsim 2$, as well as for the $ \omega^3 $ observable which loses about 1 natural time unit of predictability horizon for the 0.6 normalized RMSE threshold.  Reducing the dataset size by another order of magnitude to $n = 640$ (using a sampling interval of $ \Delta t = 0.1 $ and $\ell = 200$ eigenfunctions), the reduction of skill is, as might be expected, more noticeable, bringing down the predictability horizon for $\omega^1$ and $\omega^2$ to $\tau \simeq 1.1$ natural time units. Still, despite this reduction of skill, the ability to control the complexity of the forecast function by controlling the number of eigenfunctions allows KAF to behave stably in sparsely sampled environments.  

In Figure~\ref{figL63RMSE}, we also show RMSE results for a larger dataset consisting of $n = \text{512,000}$ samples taken at a $\Delta t = 0.01$ sampling interval. Using $\ell = 6000$ eigenfunctions, this larger dataset is seen to provide a noticeable improvement of skill over the $n =\text{64,000}$ benchmark, particularly for $ \tau \gtrsim 1$ leads. For instance, the $n = \text{512,000}$ setup maintains lower than 0.6 normalized RMSE values out to $\tau \simeq 2.7$, which represents a $\simeq 30\%$ increase over the $n=\text{64,000}$ case.

\begin{figure*}
    \includegraphics[width=\linewidth]{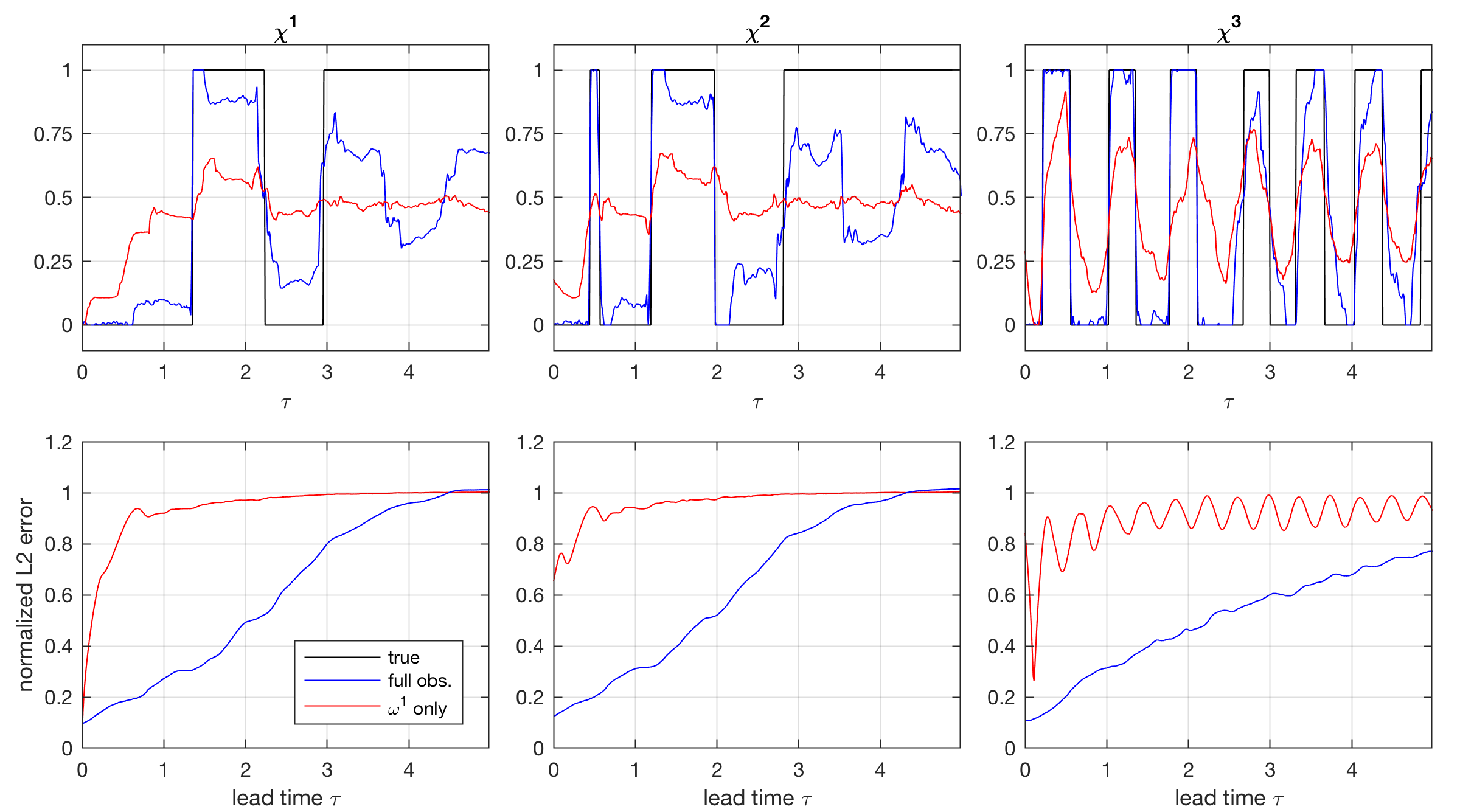}
    \caption{\label{figL63Prob}As in Figure~\ref{figL63Pred}, but for prediction of the characteristic functions $ \chi^1$, $\chi^2$, and $\chi^3 $, representing the L63 states whose components $\omega^1$, $\omega^2$, and $\omega^3$, respectively, are greater than their mean. The initial condition in the top panels is the same as in the top panels of Figure~\ref{figL63Pred}.} 
\end{figure*}

\paragraph{Forecasting conditional probability} As our final numerical results, we show in Figure~\ref{figL63Prob} trajectory and RMSE results for prediction of the characteristic functions $ \chi^i $, corresponding to the conditional probabilities for $ \omega^i $ to take greater-than-average values. These forecasts were obtained using the empirical target functions $ g_{\tau,\ell,n}$ from~\eqref{eqKAFProb}, constructed using the same parameter values as the state vector forecasts in Figure~\ref{figL63Pred} based on $f_{\tau,\ell,n}$. Compared to forecasts of the state vector components, forecasts of characteristic functions are hampered by the fact that  a characteristic function has discontinuities (apart from trivial cases), potentially inducing Gibbs oscillations in approximations by finite linear combinations of RKHS functions. Such oscillations may in turn induce overshoots outside the interval $ [ 0, 1 ] $, necessitating the use of thresholding in~\eqref{eqKAFProb}. Evidence of active thresholding can be seen in Figure~\ref{figL63Prob}, particularly at early times ($\tau \lesssim 0.5$) where the relative RMSE is significantly larger than the corresponding state vector results in Figure~\ref{figL63Pred}. Despite that, the conditional probability results are broadly consistent with their state vector counterparts. That is, $\chi^1$/ $\chi^2 $ and $ \chi^3 $  remain predictable out to 2--3 and 4--5 natural time units, similarly to  $ \omega^1$/ $\omega^2$ and $ \omega^3 $, respectively, and as expected, the fully observed forecasts fare substantially better than the partially observed ones. It is worthwhile noting that unlike $ \omega^i $, the $\chi^i$ forecasts have not converged to near-constant values at the end of the examined prediction intervals (i.e., at $ \tau = 5$).   

\section{\label{secConclusions}Conclusions}
As shown throughout this paper, the theoretical underpinning of the kernel approach to forecasting is that it approximates the conditional expectation of observables, in the sense of minimizing mean square forecast error.  The extent to which approximating the conditional expectation is one of the better ways of producing forecasts depends on the specific dynamical system and the goals of the forecaster, but is a natural and common target in many applications. Although there are many ways of achieving this approximation, we have shown in this article that the kernel approach is a distinguished such method both theoretically, given the central importance of Hilbert space theory for both kernels and the conditional expectation, and practically, as the computation requires only eigenfunction computation and matrix multiplication.

In contrast to the usual expositions of kernel methods that present kernels and RKHSs as an axiomatic starting point, we have followed a different order in which the appearance of kernels arises naturally within a learning framework (with hypothesis spaces as a prominent object) and with error minimization as a starting point. It is this perspective on kernel based forecasting, i.e. one that deemphasizes kernels in favor of conditional expectation, that is the main contribution of this paper. Additionally, we presented and proved the results that the empirically obtained kernel forecasting function approaches, in the limit of large data, the ideal kernel forecasting function, which itself approaches the true conditional expectation as more principal components (kernel eigenfunctions) are utilized. Precise estimates on the rate of convergence is an active area of research and one that depends on specific aspects of both the dynamical system and the employed kernels.

One of the advantages of an understanding of kernel forecasting based on conditional expectation, so we have argued, is that it enables the computation of a host of related quantities, including conditional probability and conditional variance. The former can be used to handle the binary classification problem that arises when trying to detect extreme or rare events. The latter, meanwhile, is instrumental in providing more informative forecasts that detail the level of uncertainty involved. Another benefit of the statistical learning framework is that it shows the connection between the two most common kernel methods, KPCR (of which KAF is an example) and KRR; in particular, they both follow the same variational logic, but the former is based on a choice of a linear hypothesis space whereas the latter uses a nonlinear one. Although KRR may be a simpler algorithm to implement, and may be more accurate in the presence of noise, KPCR can converge much more rapidly when the predictand happens to lie in the space spanned by the leading principal components. We have also shown how KAF can be implemented using a class of non-symmetric kernels satisfying a detailed-balance condition, as well as kernels based on delay-coordinate maps.

Applications of KAF to two low-dimensional dynamical systems were presented for the sake of illustration. The first system, periodic flow on a circle, is in fact not a system for which conditional expectation is a good estimate of forecasts, at least when the abscissa is the only quantity on which the forecast is conditioned. Nevertheless, we demonstrate that KAF implemented with a radial Gaussian kernel converges to this conditional expectation quite rapidly, and that the dependence of its error on the number of principal components follows a U-shaped curve that is characteristic of the classic bias-variance tradeoff of statistical learning. The second system, the L63 system, exhibits a number of the hallmark challenges in forecasting of complex systems, including invariant measures supported on complicated sets (fractal attractors) and mixing dynamics. Despite these challenges, we saw that KAF, implemented with a judiciously chosen variable-bandwidth, Markov-normalized Gaussian kernel, successfully predicts the state vector components, as well as their associated conditional probabilities to take greater-than average values. As expected, conditioning  on the full state produces better forecasts than conditioning on just partial observations of the state, but in both cases the method yielded adequate uncertainty quantification through estimates of the conditional variance. It was also found that incorporating a sufficient number of delays leads to recovery of most of the forecast skill lost due to partial observations. The L63 example also demonstrates that forecasts based on partial conditioning are better for some choices of response variables than others. In particular, as expected from symmetry considerations, the first coordinate of the state vector has greater predictive value for the third coordinate than for the second coordinate. 

There are two chief challenges in utilizing kernel methods in real-world applications. The first is an appropriate choice of response, as well as a covariate variable with sufficiently rich predictive value. The second challenge, particularly when dealing with very high-dimensional covariate spaces, is a choice of kernel such that as much of the dynamical features of interest can be characterized by as few of the leading principal components as possible. In general, the response and covariate are selected with the certainty that there is close association between the two, but with the precise nature of the correspondence being either unknown, or intractable to reproduce analytically or numerically. In real-world applications, this issue is further compounded by the fact that the response space is oftentimes multi-dimensional. While in this paper we did not directly address this situation, it is natural to consider extensions of KAF to the setting of vector-valued response functions using operator-valued kernel techniques \cite{MicchelliPontil05,SlawinskaEtAl18} for multi-task learning. Another potential direction for future research is to establish connections between KAF an RKHS embeddings of probability distributions \cite{SongEtAl09,SriperumbudurEtAl11,KlusEtAl19}. As for the choice of kernel, recent approaches for learning kernels targeted to specific response functions \cite{OwhadiYoo19} could potentially provide effective ways of ensuring that the response is well-captured by the leading eigenspaces of the corresponding integral operator, thus improving forecast skill. The main goal of this paper has been to clarify the theoretical justification for utilizing kernels in forecasting observables of dynamical systems, so that the forecaster can focus on the remaining problem of leveraging specific scientific knowledge of the system into optimal choices of response, covariate, and kernel.

\section*{Acknowledgments} Dimitrios Giannakis acknowledges support by  ONR YIP grant N00014-16-1-2649, NSF grants DMS-1521775 and 1842538, and DARPA grant HR0011-16-C-0116. Romeo Alexander was supported as a PhD student from the first NSF grant and the DARPA grant. The authors would like to thank Suddhasattwa Das, Krithika Manohar, and Andrew Stuart for fruitful conversations. In addition, they would like to thank two anonymous Reviewers for constructive comments which have led to improvements of the manuscript. Dimitrios Giannakis is grateful to the Department of Computing and Mathematical Sciences at the California Institute of Technology for hospitality and for providing a stimulating environment during a sabbatical, where a portion of this work was completed. 

\section*{Declaration of interest}

None.

\appendix

\section{\label{app}Definitions and technical results}

\subsection{\label{appPInv}Moore-Penrose pseudoinverse}

We state below the definition of the pseudoinverse of a linear map between Hilbert spaces \cite{BeutlerRoot73}. 

\begin{defn}[Moore-Penrose pseudoinverse] \label{defnPInv} Let $H_1$ and $H_2$ be Hilbert spaces over the complex numbers, and $ A : D(A) \to H_2$ a closed linear map with dense domain $D(A) \subseteq H_1$. Then, a densely defined operator $A^+ : D(A^+) \to H_1$ with domain $ D(A^+)\subseteq H_2$ is said to be a Moore-Penrose pseudoinverse of $A$ if (i) $ \ker A^+ = \ran A^\perp $; (ii) $ \overline{ \ran A^+} = \ker A^\perp $; and (iii) $ A A^+ f = f $ for all $f \in \ran A$.
\end{defn}

If $A^+$ in the above definition exists, then it is closed and unique. Moreover, if $A$ has closed range, then $A^+ $ always exists, and is a bounded operator with $D(A^+) = H_2.$ If $A$ is bounded, we can express $A^+ $ on the potentially restricted domain $D((AA^*)^+) \subseteq D(A^+)$ of the pseudoinverse of the self-adjoint operator $A A^*$ through a formula with a direct counterpart in finite-dimensional linear algebra, viz.
\begin{equation}
    \label{eqPInv}
    A^+ f = A^* ( A A^*)^+ f, \quad \forall f \in D( (AA^*)^+).
\end{equation}

\subsection{\label{appAsym}Proof of Proposition~\ref{propAsym}}

First, the symmetry and positive-definiteness of $k$ follow directly from its definition and the detailed-balance condition in~\eqref{eqBalance}. In particular, since $ d $ is strictly positive, \eqref{eqKHat} implies that $k$ is positive-definite if and only if $\hat k $ is positive-definite, and the latter is indeed positive-definite since it is related to the positive-definite kernel $w $ by a similarity transformation.

Next, the integral operator $W$ is equal to $\iota^* M_d$, and because $M_d$ is a bounded, invertible operator with bounded inverse, $W$ is well-defined on $L^2(\rho)$ and $ \ran W = \ran \iota^* \subseteq \mathcal K $, proving Claim~(i).

Turning to Claim~(ii), the fact that $J$ is a trace-class non-negative operator with real eigenvalues follows from its relation to $\hat G$ (which has all of these properties by positive-definiteness and continuity of $\hat k$) via the similarity transformation in~\eqref{eqJGHat}. In addition, the existence of the biorthonormal Riesz bases $ \{ \xi_1, \xi_2, \ldots \} $ and $ \{ \xi'_1, \xi'_2, \ldots \} $ follows from~\eqref{eqXi} in conjunction with continuity of $d^{1/2}$ and $d^{-1/2}$ (and thus boundedness of these functions on the compact support of $\rho$).

Finally, to prove Claim~(iii) note first that by Definition~\ref{defnPInv},
\begin{align*}
    D(J^+) &= \ran J \oplus \ker J^* = \ran( \iota W ) + \ker( W^* \iota^* ) \\
    &= \ran( \iota \iota^* M_d ) + \ker( M_d \iota \iota^* ) = \ran(\iota \iota^*) \oplus \ker( \iota \iota^* ) \\
    &= \ran G \oplus \ker G = D( G^+ ),
\end{align*}
where the second equality in the second line follows from the fact that $M_d$ is a bounded, invertible operator with bounded inverse. Moreover, 
\begin{displaymath}
    \ran G = \ran(\iota \iota^*) \subseteq \ran \iota,
\end{displaymath}
and because $\ran \iota^* \subseteq \ker \iota^\perp $, we have
\begin{equation}
    \label{eqRanG}
    \overline{\ran G } = \ker G = \ker \iota^* = \overline{\ran \iota}.
\end{equation}
It therefore follows that $ \ran G $ is a dense subspace of $ \ran \iota $.  Now, by Lemma~\ref{lemNystrom} and~\eqref{eqRanG}, 
\begin{displaymath}
    D(\tilde T^+) = \ran \iota \oplus \ker \iota^* = \ran \iota \oplus \ker G, 
\end{displaymath}
so we conclude that $D(J^+) = D(G^+)$  is a dense subspace of $D(\tilde T)$, as claimed. Moreover,
\begin{align*}
    \tilde T\rvert_{D(J^+)} &= \tilde T\rvert_{D(G^+)} = \iota^* G^+ \\
    &= \iota^* ( J M_d^{-1} )^+ = \iota^* M_d J^+ = W J^+.
\end{align*}
The expression for $\tilde T f $ in~\eqref{eqNystAlt} follows from the result just proved, the definitions of $\xi_i$, $\xi'_i$, and $\theta_i$ in~\eqref{eqXi} and~\eqref{eqTheta}, and the fact that $\tilde T $ is a closed operator. We then verify that the $\vartheta_i$ are indeed orthogonal, viz.
\begin{align*}
    \langle \vartheta_i, \vartheta_j \rangle_{\mathcal K} &= \frac{1}{\sqrt{\eta_i \eta_j}} \langle W M_d^{-1/2} \hat \phi_i, W M_d^{-1/2} \hat \phi_j \rangle_{\mathcal K} \\
    &=\frac{1}{\sqrt{\eta_i \eta_j}} \langle \iota^* M_d^{1/2} \hat \phi_i, \iota^* M_d^{1/2} \hat \phi_j \rangle_{\mathcal K} \\
    &=\frac{1}{\sqrt{\eta_i \eta_j}} \langle \hat\phi_i, M_d^{1/2} \iota \iota^* M_d^{1/2} \hat \phi_j \rangle_{L^2(\rho)} \\
    &=\frac{1}{\sqrt{\eta_i \eta_j}} \langle \hat \phi_i, \hat G_\rho \hat \phi_j \rangle_{L^2(\rho)} = \delta_{ij}.
\end{align*}
This completes the proof of Claim~(iii) and Proposition~\ref{propAsym}.

%\bibliography{cKAF_dg,bibliography_dg}

\end{document}